\documentclass[11pt,a4paper]{amsart}
\usepackage{amsmath,amssymb,bbm,epsfig}
\usepackage[all]{xy}

\newtheorem{theorem}{Theorem}[section] 
\newtheorem{lemma}[theorem]{Lemma}
\newtheorem{proposition}[theorem]{Proposition}
\newtheorem{propdef}[theorem]{Proposition-Definition}
\newtheorem{corollary}[theorem]{Corollary}
\newtheorem{definition}[theorem]{Definition}
\newtheorem{remark}[theorem]{Remark}

\newtheorem{rmq}[theorem]{Remark}

\newenvironment{pr}{\noindent {\bf Proof. }}{\ \ \ \hfill $\square$ \\}

\newcommand{\bxp}{\boxplus}

\def\vfi{\varphi}
\def\pro{probability }
\def\be{\begin{equation}}
\def\ee{\end{equation}}
\def\bea{\begin{eqnarray}}
\def\eea{\end{eqnarray}}
\def\beas{\begin{eqnarray*}}
\def\eeas{\end{eqnarray*}}

\def\lf{\left}
\def\ri{\right}
\def\Tr{\operatorname{Tr}}
\def\tr{\operatorname{tr}}

\newcommand{\lan}{\langle}
\newcommand{\ran}{\rangle}
\newcommand{\f}{\frac}
\newcommand{\ff}{\frac{1}}
\newcommand{\ste}{\, ;\, }
\newcommand{\st}{such that }
\newcommand{\eps}{\varepsilon}

\newcommand{\ud}{\mathrm{d}}

\newcommand{\mc}{\mathcal}
\def\B{{\mathcal B}}

\def\M{{\mathcal M}}

\def\A{{\mathcal A}}
\def\C{{\mathbb C}}
\def\E{{\mathbb E}}
\def\Herm{{\mathcal H}}

\def\R{{\mathbb R}}
\def\Sym{{\mathfrak S}}

\def\U{{\mathbb U}}

\title[From classical to free independence]{A continuous semigroup of notions of independence between the classical and the free one}
\author{Florent Benaych-Georges and Thierry L\'evy}
\keywords{Free Probability ; Independence ; Random Matrices ; Unitary Brownian Motion ; Convolution ; Cumulants}
\subjclass[2000]{46L54, 15A52}
\address{Florent Benaych-Georges, LPMA, UPMC Univ Paris 6, Case courier 188, 4, Place Jussieu, 75252
Paris Cedex 05, France, and CMAP, ´ Ecole Polytechnique, route de Saclay, 91128 Palaiseau Cedex,
France}
\email{florent.benaych@gmail.com}
\address{Thierry L\'evy, CNRS and \'Ecole Normale Sup\'erieure, DMA, 45, rue d'Ulm, F-75005 Paris}
\email{levy@dma.ens.fr}
\begin{document}
\maketitle

\begin{abstract} In this paper, we investigate a continuous family of notions of independence which interpolates between the classical and free ones for non-commutative random variables. These notions are related to the liberation process introduced by D. Voiculescu. To each notion of independence correspond new convolutions of probability measures, for which we establish formulae and of which we compute simple examples. We prove that there exists no reasonable analogue of classical and free cumulants associated to these notions of independence. 
\end{abstract}

\tableofcontents

\section*{Introduction}

Let $\mu$ and $\nu$ be two Borel probability measures on the real line $\R$. The classical convolution of $\mu$ and $\nu$ is the probability measure on $\R$, denoted by $\mu *\nu$, which is the distribution of the sum of two classical independent random variables with respective distributions $\mu$ and $\nu$. Let us describe $\mu*\nu$ in an alternative way. To each $n\times n$ matrix with eigenvalues $\lambda_{1},\ldots,\lambda_{n}$, we associate its spectral measure, which is the probability measure $\frac{1}{n}\sum_{i=1}^{n}\delta_{\lambda_{i}}$. Let $(A_{n})_{n\geq 1}$ and $(B_{n})_{n\geq 1}$ be two sequences of diagonal real matrices, with $A_{n}$ and $B_{n}$ of size $n$ for all $n\geq 1$, such that the spectral measure of $A_{n}$ (resp. of $B_{n}$) converges, as $n$ tends to infinity, to $\mu$ (resp. to $\nu$). For each $n\geq 1$, let $S_{n}$ be a random matrix chosen uniformly among the $n!$ permutation matrices of size $n$. Then the spectral measure of $A_{n}+S_{n}B_{n}S_{n}^{-1}$ converges, as $n$ tends to infinity, to $\mu*\nu$.

If we replace, for each $n\geq 1$,  the matrix $S_{n}$ by a random matrix $U_{n}$ chosen in the unitary group $U(n)$ according to the Haar measure, then the spectral measure of $A_{n}+U_{n}B_{n}U_{n}^{-1}$ converges, as $n$ tends to infinity, to the free convolution of $\mu$ and $\nu$, a probability measure on $\R$ denoted by $\mu \boxplus \nu$.

This way of describing classical and free convolutions suggests a natural way to interpolate between them. Indeed, consider, for all $n\geq 1$, a properly scaled Brownian motion $(U_{n,t})_{t\geq 0}$ issued from the identity matrix on the unitary group $U(n)$. Given $t\in [0,+\infty)$, one may consider the spectral measure of $A_{n}+U_{n,t}S_{n} B_{n} S_{n}^{-1}U_{n,t}^{-1}$, and ask for the limit of this distribution as $n$ tends to infinity. For $t=0$, the matrix $U_{n,0}$ is the identity matrix and we find the classical convolution of $\mu$ and $\nu$. For $t=+\infty$, that is, when $U_{n,t}$ is replaced by its limit in distribution as $t$ tends to infinity, which is a uniformly distributed unitary matrix, we recover the free convolution of $\mu$ and $\nu$. For any other $t\in (0,+\infty)$, it turns out that one finds a probability measure which depends only on $\mu$, $\nu$ and $t$ and which we denote by $\mu *_t\nu$. Note that this definition of $*_{t}$ can be considered as a particular case of the so-called {\em liberation process} introduced by Voiculescu \cite{v99}.

Consider for example the case where $\mu=\nu=\frac{1}{2}(\delta_1+\delta_{-1})$. Then $\mu*\nu=\frac{1}{4}\delta_{-2}+\frac{1}{2} \delta_0 + \frac{1}{4} \delta_2$ and it is well known that $\mu\boxplus \nu={\mathbbm 1}_{[-2,2]}(x) \frac{dx}{\pi \sqrt{4-x^2}}$, a dilation of the arcsine law \cite[Example 12.8]{NS}. One may wonder which probability measures interpolate between $\mu*\nu$ and $\mu\boxplus \nu$. We will prove that
$$\forall t\geq 0\; , \; \; \frac{\delta_1+\delta_{-1}}{2} *_{t} \frac{\delta_1+\delta_{-1}}{2} = {\mathbbm 1}_{[-2,2]}(x)  \frac{\rho_{4t}(e^{4i \arccos \frac{x}{2}})}{\pi \sqrt{4-x^2}} \; dx.$$
Here, for all $t> 0$ and $\theta\in \R$, $\rho_{t}(e^{i\theta})$ is the density at $e^{i\theta}$, with respect to the uniform probability measure on the unit circle, of the distribution of the free unitary Brownian motion at time $t$. This distribution is also the limit, as $n$ tends to infinity, of the spectral measure of $U_{n,t}$. There is no simple formula for this distribution, which apparently has to be taken as a fundamental function in any problem involving the large $n$ asymptotics of the Brownian motion on the unitary group $U(n)$. However the moments of this distribution are known since P. Biane first computed them \cite{b97}. It follows for instance from the previous expression that $\mu*_t \nu$ has a density with respect to the Lebesgue measure for all $t>0$ and that its support, which one can compute for all $t\geq 0$, is the whole interval $[-2,2]$ if and only if $t\geq 1$. \\

The family of operations $*_{t}$ is really just a by-product of a more fundamental  construction, which is that of a continuous family of independence (or dependence) structures between non-commutative random variables which interpolates between classical independence and freeness. Indeed, we will define, for all $t\in [0,+\infty]$, a notion of independence between two subalgebras of a non-commutative probability space, which we call {\em $t$-freeness} and which, for $t=0$ (resp. $t=+\infty$), coincides with classical independence (resp. freeness). Once this structure is defined, it is straightforward to define additive or multiplicative convolution of $t$-free self-adjoint or unitary elements, thus giving rise to several  operations on probability measures: additive or multiplicative convolution of probability measures with compact support on  $\R$, denoted by $*_{t}$ and $\odot_{t}$ ; multiplicative convolution of probability measures on the unit circle, also denoted by $\odot_{t}$.

The idea of seeking a continuous way of passing from classical to free independence is presumably as old as the theory of free probability itself, but the research of such a continuum has been broken off by  a paper of Roland Speicher in 1997 \cite{s97}, where he has shown that no other notion of independence than the classical and the free ones can be the base of a {\em reasonable} probability theory. Indeed $t$-freeness does not satisfy all the axioms enforced by R. Speicher because it is not an {\em associative} notion of independence. This axiom of associativity states, roughly, that if $X,Y, Z$ are three random variables such that $X$ is independent of $Y$ and $Z$ is independent of $\{X,Y\}$, then $X$ must be independent of $\{Y,Z\}$. Instead of this, what is true with $t$- freeness is that for all $s,t\geq 0$, if  $X,Y, Z$ are three random variables such that $X$ is $t$-free with $Y$ and $Y$ is $s$-free with $Z$, then under certain additional hypotheses, $X$ will be $(s+t)$-free with $Z$. This is of course related to the semi-group property of the Brownian motion. \\

There are several ways to characterize and deal with independence and freeness. The first one, which we have already mentioned, is to relate them with matrix models. The second one is to describe them by means of computation rules: the expectation factorizes with respect to independent subfamilies of random variables, whereas the expectation of a product of free elements can be computed using the fact that if $x_1,\ldots, x_n$ are centered and successively free, then their product is centered. The third way to describe independence and freeness is to identify integral transforms which linearize them (namely the logarithm of Fourier transform or the $R$-transform). This amounts to describing classical and free cumulants. The last way, a bit more abstract, is to consider tensor or free products: a family of random variables is independent (resp. free) if and only if it can be realized on a tensor product (resp. free product) of probability spaces. 
 
In the present paper, we look for the analogues of all these approaches for the notion of $t$-freeness. We begin, in Section \ref{10.9.08.1}, by giving the definition of a $t$-free product and presenting  the corresponding random matrix model. Then, in Section \ref{10.9.08.2}, we state the computation rules, which are best understood as a family of differential equations. Finally, in Section \ref{10.9.08.3}, we prove that no notion of cumulants of order greater than $6$ can be associated to the notion of $t$-freeness. More precisely, we show that there does not exist a universally defined $7$-linear form on any non-commutative probability space with the property that this form vanishes whenever it is evaluated on arguments which can be split into two non-empty subfamilies which are $t$-free, unless $t=0$ or $t=+\infty$. This can be summarized in the following diagram.  \\

\begin{figure}[h!]
\begin{center}
\begin{tabular}{c|c|c|c|c|}
\cline{2-5}
 &  Matrix model & Computation rules & Cumulants & Algebraic structures \cr
\hline
\multicolumn{1}{|c|}{Indep.} & $A+SBS^{-1}$ & Factorization & Class. cumulants & Tensor product \cr
\hline
\multicolumn{1}{|c|}{$t$-freeness} & $A+U_tS B S^{-1}U_t^{-1}$ & Differential system & {\bf Do not exist} & $t$-free product \cr
\hline
\multicolumn{1}{|c|}{Freeness} & $A+UBU^{-1}$ & $\varphi(x_1\ldots x_n)=0$ & Free cumulants & Free product \cr
\hline
\end{tabular}
\caption{The main computation rule for freeness is that $\varphi(x_1\ldots x_n)=0$ as soon as $x_1,\ldots,x_n$ are successively free and centered. For $t$-freeness, the computation rules are best expressed as a differential system relating the distributions of the $t$-free products of two families of random variables for different values of $t$.}
\end{center}
\end{figure}

\section{Preliminaries}

In this section, we review the notions of non-commutative probability which are relevant to the definition of $t$-freeness.

\subsection{Probability space, distribution}
Non-commutative probability is based on the following generalization  of the notion of probability space. 

\begin{definition}[Non-commutative \pro space]
A {\em non-commutative probability space} is  a pair $(\mc{A},\vfi)$, where: \begin{itemize} \item $\A$ is an algebra over $\C$ with a unit element denoted by $1$, endowed with an operation of adjunction  $x\mapsto x^*$ which is $\C$-antilinear, involutive and satisfies $(xy)^*=y^*x^*$ for all $x,y\in \A$, \item $\vfi\,:\, \A\to \C$ is a linear form on $\A$,  satisfying $\vfi(1)=1$, $\vfi(xy)=\vfi(yx)$, $\vfi(x^*)=\overline{\vfi(x)}$ and $\vfi(xx^*)\geq 0$ for all $x,y\in \A$.\end{itemize} 
The linear form $\vfi$ is often called the {\it expectation} of the non-commutative \pro space.
\end{definition}

Two fundamental examples are the algebra $L^{\infty-}(\Omega,\Sigma,\mathbb P)$ of complex-valued random variables with moments of all orders on a classical probability space, endowed with the complex conjugation and the expectation (we will say that this non-commutative probability space is {\em inherited} from $(\Omega,\mathcal A,\mathbb P)$); and the algebra $\M_{n}(\C)$ endowed with the matricial adjunction and the normalized trace.

%Here are two very natural examples of non-commutative probability spaces.
%\begin{example} [Classical \pro space]\label{23.06.08.11h35}\rm 
%{\rm Let $(\Omega, \Sigma, \mathbb{P})$ be a classical \pro space where the expectation is denoted by $\mathbb{E}$. Then the algebra of complex-valued random variables defined on this space having moments to all orders, endowed with the involution $X^*=\overline{X}$ and with the linear form $\mathbb{E}$ is a non-commutative \pro space, of which we will say that it is {\em inherited} from a classical probability space.} \end{example} 
 
%\begin{example}[Space of matrices endowed with the normalize trace]{\rm 
%Let $n\geq 1$ be an integer. Then the algebra $\A=\M_n(\C)$ endowed with the linear form $\vfi$ which is the normalized trace is a non-commutative probability space.}\end{example}
 
%Note that  the previous  example can be generalized by replacing $\M_n(\C)$ by a finite von Neumann algebra and the normalized trace by one of its tracial states. 
   
%The second notion from classical \pro which has a natural non-commutative generalization is the one of distribution of a family of random variables. 

\begin{definition}[Non commutative distribution]\label{def ncd}
Let $(\A,\varphi)$ be a non-commutative probability space. The {\it non-commutative distribution} of a family $(a_1, \ldots, a_n)$  of elements of $\A$ with respect to $\varphi$ is the linear map defined on the space of polynomials in the non-commutative variables $X_1, X_1^*, \ldots, X_n, X_n^*$ which maps any such polynomial $P$  to  $\vfi(P(a_1,a_1^*,\ldots, a_n,a_n^*))$.
\end{definition}

%\begin{rmq}\label{suite.moments}{\rm Consider a self-adjoint element $a$ (i.e. \st $a=a^*$) in a non-commutative \pro space $(\A, \vfi)$. Then by the hypothesis $\vfi(xx^*)\geq 0$ for all $x\in \A$, the distribution of $a$ is a linear form on $\C[X]$ which is non negative on the polynomials which are non negative on the real line. Hence it is the integration with respect to a \pro measure on the real line. Note that this measure is unique if and only if it is determined by its moments, which is the case when the measure has compact support, i.e. when there is a constant $M$ \st for all $n\geq 0$, one has $\vfi(a^{2n})\leq M^{2n}$.}
%\end{rmq}

The link between the classical notion of distribution and the non-commutative one is the following. Consider a self-adjoint element $a$ in a non-commutative \pro space $(\A, \vfi)$, that is, an element such that $a=a^*$. Since $\vfi(xx^*)\geq 0$ for all $x\in \A$, the distribution of $a$ is a linear form on $\C[X]$ which is non-negative on the polynomials which are non-negative on the real line. Hence, it can be represented as the integration with respect to a \pro measure on the real line. This probability measure is unique if and only if it is determined by its moments, which is in particular the case when it has compact support, or equivalently when there exists a constant $M$ \st for all $n\geq 0$, one has $\vfi(a^{2n})\leq M^{2n}$. 

Similarly, the distribution of a unitary element $u$, that is, an element such that $uu^{*}=u^{*}u=1$, is the integration with respect to a probability measure on the unit circle of $\C$. Since the circle is compact, there is no issue of uniqueness in this case.

%\begin{rmq}\label{suite.moments.bis}{\rm The content of the previous remark can be generalized to any unitary element: its distribution is the integration with respect to a probability measure on the unit circle of $\C$.}
%\end{rmq}

%\begin{rmq}{\rm Note, as a generalization of   Remark \ref{suite.moments}, that in the case of a non-commutative \pro space inherited from a classical one (see Example \ref{23.06.08.11h35}), the non-commutative distribution of a family of random variables is the integration of polynomials with respect to the law of the family (defined in the classical way).}
%\end{rmq}

\subsection{Independence, freeness and random matrices} \subsubsection{Definitions and basic properties}

%After the notions  of  \pro space and of  distribution, one of the most important notion in \pro theory  is the independence. 

We shall recall the definitions of the two notions of independence in a non-commutative probability space between which our main purpose is to interpolate. The first one is a straightforward translation of the classical notion of independence in the non-commutative setting, which coincides with the original notion in the case of a non-commutative \pro space inherited from a classical one. The second one is the notion of freeness, as defined by Voiculescu \cite{vdn91}, which is called freeness. 

In this paper, by a subalgebra of the algebra of a non-commutative \pro space, we shall always mean a subalgebra which contains $1$ and which is stable under the operation $x\mapsto x^*$. 

\begin{definition}[Independence and freeness]\label{23.06.08.17h04}
Let  $(\mc{M},\vfi)$ be a non-commutative \pro space. The kernel of $\vfi$ will be called the set of {\it centered elements}. Consider a family $(\A_i)_{i\in I}$ of subalgebras of $\M$. \begin{itemize}\item The family $(\A_i)_{i\in I}$ is said to be {\it independent} if \begin{itemize}\item[(i)] for all $i\neq j\in I$, $\A_i$ and $\A_j$ commute, \item[(ii)] for all $n\geq 1$, $i_1, \ldots, i_n\in I$ pairwise distinct, for all family $(a_1, \ldots, a_n)\in \A_{i_1}\times \cdots\times \A_{i_n}$ of centered elements, the product $a_1\cdots a_n$ is also centered.\end{itemize}
 \item The family $(\A_i)_{i\in I}$ is said to be {\it free} if for all $n\geq 1$, $i_1, \ldots, i_n\in I$ \st $i_1\neq i_2, i_2\neq i_3, \ldots, i_{n-1}\neq i_n$, for all family $(a_1, \ldots, a_n)\in \A_{i_1}\times \cdots\times \A_{i_n}$ of centered elements, the product $a_1\cdots a_n$ is also centered.
\end{itemize}
\end{definition}

On a classical probability space $(\Omega,\Sigma,\mathbb P)$, a family $(\Sigma_{i})_{i\in I}$ of sub-$\sigma$-fields of $\Sigma$ is independent with respect to $\mathbb P$ if and only if the subalgebras $\left(L^{\infty-}(\Omega,\Sigma_{i},\mathbb P)\right)_{i\in I}$ of $\left(L^{\infty-}(\Omega,\Sigma,\mathbb P),\mathbb E\right)$ are independent in the sense of the definition above.

%Note that in the case  of a non-commutative \pro space inherited from a classical one (see Example \ref{23.06.08.11h35}),  the notion of independence as it is presented here is exactly the classical notion of independence: with the notation of Example \ref{23.06.08.11h35}, let us define, for all family $(\Sigma_i)_{i\in I}$ of sub-$\sigma$-fields of $\Sigma$, and for all $i\in I$, the set $\A_i$ of random variables of $\A$ which are $\Sigma_i$-measurable. Then the family $(\Sigma_i)_{i\in I}$ is independent in the classical sense if and only if the family $(\A_i)_{i\in I}$ is independent in the sense of Definition \ref{23.06.08.17h04}. 

In the classical setting again, a family of random variables is independent if and only if its joint distribution is the tensor product of the individual ones. In the following definition and proposition, we translate this statement into our vocabulary, and give its analogue for freeness. These definitions prepare those which we will give later for $t$-freeness. 
 
\begin{definition}[Tensor and free product]\label{23.06.08.17h09}{ Let  $(\mc{A}_1,\vfi_1)$ and $(\mc{A}_2,\vfi_2)$ be two  non-commutative \pro spaces. \begin{itemize}\item Their  {\em tensor product}, denoted by $(\mc{A}_1,\vfi_1)\otimes(\mc{A}_2,\vfi_2)$, is the non-commutative \pro space with algebra the tensor product of  unital algebras $\mc{A}_1\otimes\mc{A}_2$,  on which the adjoint operation and the expectation are   defined by  $$\forall (x_1,x_2)\in \A_1\times \A_2, (x_1\otimes x_2)^*=x_1^*\otimes  x_2^*,\quad\vfi(x_1\otimes x_2)=\vfi_1(x_1)\vfi_2(x_2).$$\item Their  {\em free product}, denoted by  $(\mc{A}_1,\vfi_1)*(\mc{A}_2,\vfi_2)$, is the non-commutative \pro space with algebra the free product of  unital algebras $\mc{A}_1*\mc{A}_2$, with adjoint operation and expectation  defined uniquely by the fact that for all $n\geq 1$, for all $i_1\neq \cdots\neq i_n\in \{1,2\}$, for all $(x_1, \ldots, x_n)\in \A_{i_1}\times\cdots\times\A_{i_n}$, $$(x_1\cdots x_n)^*=x_n^*\cdots x_1^*$$ and $x_1\cdots x_n$ is centered whenever all $x_i$'s are.
\end{itemize}}\end{definition}

This definition can easily be extended to products of finite or infinite families of non-commutative \pro spaces, but we have restricted ourselves to what is needed in this article. We can now explain the link between these products and the notions of independence and freeness.

\begin{proposition}[Characterization of independence and freeness]\label{23.06.08.18h29}  Let  $(\mc{M},\vfi)$ be a    non-commutative \pro space. Let $\A_1,\A_2$ be subalgebras of $\A$.  Then the family $(\A_1,\A_2)$ is 
\begin{itemize}
\item independent if and only if $\A_1$ commutes with $\A_2$ and the unique algebra morphism defined from $\A_1\otimes \A_2$ to $\mc{M}$  which, for all $(a_1,a_2)\in \A_1\times \A_2$, maps  $a_1\otimes 1$ to $a_1$ and $1\otimes a_2$ to $a_2$, preserves the expectation from $(\A_1,\vfi_{|\A_1})\otimes  (\A_2,\vfi_{|\A_2})$ to $(\mc{M},\vfi)$,
\item  free if and only if  the unique algebra morphism defined from the free product of unital algebras  $\A_1* \A_2$ to $\mc{M}$ which, restricted to $\A_1\cup\A_2$ is the canonical injection, preserves the expectation from $(\A_1,\vfi_{|\A_1})*  (\A_2,\vfi_{|\A_2})$ to $(\mc{M},\vfi)$.
\end{itemize}
\end{proposition}  

Let us finally recall the definition of the free analogue of the classical convolution, which is meaningful thanks to the last proposition.

\begin{definition}[Additive free convolution]\label{23.06.08.20h09}{Let $\mu$ and $\nu$ be two probability measures on $\R$. The distribution of the sum of two free self-adjoint elements with respective distributions $\mu$ and $\nu$ depends only on $\mu$ and $\nu$ and will be called the {\rm free additive convolution} of $\mu$ and $\nu$, and be denoted by $\mu\bxp\nu$. 
}\end{definition}

\subsubsection{Asymptotic behavior of random matrices}
In this section, we recall matrix models for the classical and free convolution. The main notion of convergence which is involved is the following.

\begin{definition}[Convergence in non-commutative distribution]\label{23.06.08.17h10}{ Let $p$ be a positive integer and let, for each $n\geq 1$, $(M(1,n), \ldots, M(p,n))$ be a family of $n\times n$ random matrices. This family is said to {\em converge in non-commutative distribution} if its non-commutative distribution converges in \pro to a non random one, that is, if the normalized trace of any word in the $M(i,n)$'s and the $M(i,n)^*$'s converges in \pro to a constant.}\end{definition}

%As mentioned above,  random variables which are independent in the usual sense are elements of independent algebras in the sense developed in this paper. It is not so easy to construct free random variables without using explicitly the free product of definition \ref{23.06.08.17h09}. However,  in 91, Voiculescu, in \cite{v91}, proved that random matrices conjugated by independent unitary random matrices with Haar distribution are asymptotically free. 
%In the following theorem, we present this result and its analogue for freeness, where unitary matrices are replaced by permutation matrices.

\begin{theorem}[Asymptotic independence and asymptotic freeness]\label{24.06.08.1} Let us fix $p, q\geq 1$. For each $n\geq 1$, let $\mc{F}_n=(A(1,n), \ldots, A(p,n),B(1,n),\ldots, B(q,n))$ be a family of $n\times n$ random
matrices and assume that the sequence $(\mathcal F_n)_{n\geq 1}$ converges in non-commutative distribution. Assume also that for all $r\geq 1$, the entries of these random matrices are  uniformly bounded in $L^r$.
\begin{itemize}
\item Assume that these matrices are diagonal and consider, for each $n$,  the matrix $S_n$ of a uniformly distributed random permutation of $\{1,\ldots, n\}$ independent of the family $\mc{F}_n$. Then the family \be\label{26.11.08.11h40}(A(1,n), \ldots, A(p,n),S_nB(1,n)S_n^{-1},\ldots, S_nB(q,n)S_n^{-1})\ee converges in distribution to the distribution of a commutative family $(a_1,\ldots, a_p, b_1, \ldots, b_q)$ of elements of a non-commutative \pro space \st the algebras generated by $\{a_1, \ldots, a_p\}$ and $\{b_1, \ldots, b_q\}$ are independent. 
\item Consider, for each $n$,  the matrix $U_n$ of a uniformly distributed  random unitary $n$ by $n$ matrix  independent of the family $\mc{F}_n$. Then the family 
\begin{equation}
(A(1,n), \ldots, A(p,n),U_nB(1,n)U_n^{-1},\ldots, U_nB(q,n)U_n^{-1})
\end{equation}
converges in distribution to the distribution of a family $(a_1,\ldots, a_p, b_1, \ldots, b_q)$ of elements of a non-commutative \pro space \st the algebras generated by $\{a_1, \ldots, a_p\}$ and $\{b_1, \ldots, b_q\}$ are free. 
\end{itemize}
\end{theorem}

\begin{rmq}{\rm The hypothesis of uniform boundedness of the entries of the matrices in each $L^r$ could be sharply weakened for the first part of the theorem if,  instead of asking for the convergence of the non-commutative distribution of the family \eqref{26.11.08.11h40}, one would ask for the weak convergence of  the empirical joint spectral measure. This would amount to choosing, as set of test functions, the set of bounded continuous functions of $p+q$ variables instead of the set of polynomials in $p+q$  variables   (see \cite{bgchapon}, where this is precisely proved).}\end{rmq}

The first part of this theorem is much simpler than the second but seems to be also less well-known. It is in any case harder to locate a proof in the literature, so that we offer one. We shall need the following lemma. We denote by $\Vert\cdot \Vert_{2}$ the usual Hermitian norm on $\C^{n}$.

\begin{lemma}Let, for each $n\geq 1$, $x(n)=(x_{n,1},\ldots, x_{n,n})$ and $y(n)=(y_{n,1},\ldots, y_{n,n})$ be two complex random vectors defined on the same \pro space \st the random variables $$\overline{x(n)}=\frac{x_{n,1}+\cdots + x_{n,n}}{n},\quad\overline{y(n)}=\frac{y_{n,1}+\cdots+ y_{n,n}}{n}$$ converge in \pro to constant limits $x,y$ as $n$ tends to infinity. Suppose moreover that the sequences $\ff{n}\Vert x(n)\Vert_2^2$ and $\ff{n}\Vert y(n)\Vert_2^2$ are bounded in $L^2$.   Consider, for all $n$,  a uniformly distributed random permutation $\sigma_n$  of $\{1,\ldots, n\}$, independent of $(x(n), y(n))$, and define $y_{\sigma_n}(n):=(y_{n,\sigma_n(1)},\ldots,y_{n,\sigma_n(n)})$.  Then  the scalar product  
$$\ff{n}\lan x(n), y_{\sigma_n}(n)\ran=\frac{x_{n,1}y_{n,\sigma_n(1)}+\cdots+ x_{n,n}y_{n,\sigma_n(n)}}{n}$$
converges in \pro to $xy$ as $n$ tends to infinity.
\end{lemma}

\begin{pr} First of all, note that  one can suppose that for all $n$, $ \overline{x(n)}=\overline{y(n)}=0$ almost surely. Indeed, if the result is proved under this additional hypothesis, then since for all $n$, one has
$$\ff{n}\lan x(n), y_{\sigma_n}(n)\ran=\ff{n}\lan x(n)-\overline{x(n)}\cdot1_n, y_{\sigma_n}(n)-\overline{y(n)}\cdot1_n\ran+\overline{x(n)}\cdot\overline{y(n)},\quad\textrm{(with $1_n=(1,\ldots, 1)$)},$$
the result holds for general $x(n), y(n)$. So we henceforth assume that  for all $n$, $ \overline{x(n)}=\overline{y(n)}=0$. The equality $\overline{y(n)}=0$ implies, for all $n$ and all $i,j=1,\ldots, n$, that 
$$\E[y_{n, \sigma_n(i)}y_{n, \sigma_n(j)}\,|\,x(n),y(n)]=\begin{cases}\ff{n}\Vert y(n)\Vert_2^2&\textrm{if $i=j$,}\\  -\ff{n(n-1)}\Vert y(n)\Vert_2^2&\textrm{if $i\neq j$.}
\end{cases}$$
Then, using  the fact that $ \overline{x(n)}=0$, we have \beas\E\lf[\ff{n^2}\lan x(n), y_{\sigma_n}(n)\ran^2\ri]&=&\E\lf[\ff{n^3(n-1)}\Vert x(n)\Vert_2^2\Vert y(n)\Vert_2^2+\ff{n^3}\Vert x(n)\Vert_2^2\Vert y(n)\Vert_2^2\ri]\\  &=& O\lf(\ff{n}\ri),\eeas which completes the proof.
\end{pr}\\

{\noindent {\bf Proof of Theorem \ref{24.06.08.1}. }}The second point is a well-known result of Voiculescu (see \cite{vdn91}).
To prove the first one, we shall prove that the normalized trace any word in the random matrices $A(1,n), \ldots, A(p,n),S_nB(1,n)S_n^{-1},\ldots, S_nB(q,n)S_n^{-1}$ converges to a constant which is the product in two terms: the limiting normalized trace of the $A(i,n)$'s and the $A(i,n)^*$'s which appear in the word on one hand and the limiting normalized trace of the $B(j,n)$'s  and the $B(j,n)^*$'s which appear in the word on the other hand. Since the $A(i,n)$'s, the $A(i,n)^*$'s,  the $S_nB(j,n)S_n^{-1}$'s  and the $S_nB(j,n)^*S_n^{-1}$'s commute, are uniformly bounded and their non-commutative distribution converges, this amounts to proving that if $M(n), N(n)$ are two diagonal  random matrices with entries uniformly bounded in $L^r$ for all $r\geq 1$, whose normalized traces converge in \pro to constants $m,n$, then for $S_n$ the matrix of a uniform random permutation of $\{1, \ldots, n\}$ independent of $(M(n), N(n))$, the normalized trace of $M(n)S_nN(n)S_n^{-1}$ converges to $mn$.  This follows directly from the previous lemma and the proof is complete.
{\ \ \hfill $\square$}

\begin{corollary}[Matricial model for classical and free convolutions] Let $\mu,\nu$ be two \pro measures on the real line. Let, for each $n\geq 1$, $M_n,N_n$ be $n$ by $n$ diagonal random matrices with empirical spectral measures converging weakly in \pro to $\mu$ and $\nu$respectively. For each $n\geq 1$, let $S_n$ (resp. $U_n$) be a uniformly distributed $n$ by $n$ permutation (resp. unitary) random matrix independent of $(M_n,N_n)$. Then \begin{itemize}\item the empirical spectral measure of $M_n+S_nN_nS_n^{-1}$ converges weakly in \pro to the classical convolution $\mu *\nu$ of $\mu$ and $\nu$, \item the empirical spectral measure of $M_n+U_nN_nU_n^{-1}$ converges weakly in \pro to the free convolution $\mu \bxp\nu$ of $\mu$ and $\nu$.\end{itemize}
\end{corollary}

\begin{pr}In the case where $\mu,\nu$ have compact supports and the entries of the diagonal matrices $M_n,N_n$ are uniformly bounded, it is a direct consequence of the previous theorem. The general case can easily be deduced using functional calculus, like in the proof of  Theorem 3.13 of \cite{bg07}. \end{pr}

\subsection{Unitary Brownian motion, free unitary Brownian motion}

In this paragraph, we give a brief survey of the definition and the main convergence result for the
Brownian motion on the unitary group.

Let $n\geq1$ be an integer. Let $\Herm_{n}$ denote the $n^2$-dimensional real linear subspace of $\M_{n}(\C)$ which consists of Hermitian matrices. On $\M_{n}(\C)$, we denote by $\Tr$ the usual trace and by $\tr=\frac{1}{n}\Tr$ the normalized trace. Let us endow $\Herm_{n}$ with the scalar product
$\langle\cdot,\cdot\rangle$ defined by
$$\forall A,B \in \Herm_{n}\; , \;\; \langle A,B\rangle =n\Tr(A^*B)=n\Tr(AB).$$
There is a linear Brownian motion canonically attached to the Euclidean space
$(\Herm_{n},\langle\cdot,\cdot\rangle)$. It is the unique Gaussian process $H$ indexed by $\R_{+}$
with values in $\Herm_{n}$ such that for all $s,t\in\R_{+}$ and all $A,B\in \Herm_n$, one has
$$\E[\langle H_s,A\rangle \langle H_t,B\rangle]=\min(s,t) \langle A,B\rangle.$$ 

Let us consider the following stochastic differential equation:
$$U_{0}=I_{n} \; , \;\; \ud U_{t}=i (\ud H_{t})U_{t}-\frac{1}{2} U_{t} \ud t,$$
where $(U_{t})_{t\geq 0}$ is a stochastic process with values in $\M_{n}(\C)$. This linear equation admits a strong solution. The process $(U^*_{t})_{t\geq 0}$, where
$U^*_{t}$ denotes the adjoint of $U_{t}$, satisfies the stochastic differential equation
$$U^*_{0}=I_{n} \; , \;\; \ud U^*_{t}=-iU^*_{t} \ud H_{t}-\frac{1}{2} U^*_{t} \ud t.$$
An application of It\^{o}'s formula to the process $U_{t}U_{t}^*$ shows that, for all $t\geq 0$,
$U_{t}U^*_{t}=I_{n}$. This proves that the process $(U_{t})_{t\geq 0}$ takes its values in the
unitary group $U(n)$.

\begin{definition}{\rm The process $(U_{t})_{t\geq 0}$ is called the {\em unitary Brownian motion of
dimension $n$}.}
\end{definition}

As $n$ tends to infinity, the unitary Brownian motion has a limit in distribution which we now describe.
For all $t\geq 0$, the numbers
$$e^{-\frac{kt}{2}} \sum_{j=0}^{k-1} \frac{(-t)^j}{j!} \binom{k}{j+1} k^{j-1} \; , \; \; k\geq 0,$$
are the moments of a unique   probability measure on the set $\U=\{z\in\C:|z|=1\}$ invariant by the complex conjugation. We denote this probability measure by $\nu_{t}$. The following definition was given by P. Biane in \cite{b97}.

\begin{definition}\label{20.06.08.1}{ Let $(\A,\tau)$ be a  non-commutative probability space. We say that
a collection $(u_{t})_{t\geq 0}$ of unitary elements of $\A$ is a {\em free unitary Brownian
motion} if the following conditions hold.
\begin{itemize}
\item For all $s,t\geq 0$ such that $s\leq t$, the distribution of $u_{t}u_{s}^*$ is the
probability measure $\nu_{t-s}$.
\item For all positive integer $m$, for all $0\leq t_{1}\leq t_{2}\leq \ldots\leq t_{m}$, the elements
$u_{t_{1}}u_{0}^*,$ $u_{t_{2}}u_{t_{1}}^*,\ldots,u_{t_{m}}u_{t_{m-1}}^*$ are free.
\end{itemize}
}
\end{definition}

In the same paper, P. Biane has proved the following convergence result.

\begin{theorem}\label{23.06.08.1} For each $n\geq 1$, let $(U_{n,t})_{t\geq 0}$ be a Brownian motion on the unitary group $U(n)$. As $n$ tends to infinity, the collection of random matrices
$(U_{n,t})_{t\geq 0}$ converges in non-commutative  distribution to a free
unitary Brownian motion.
\end{theorem}

\section{A continuum of notions of independence}\label{10.9.08.1}

In this section, we shall define a family indexed by a real number $t\in [0,+\infty]$ of relations between two subalgebras of a  non-commutative \pro space which passes  from  the classical independence (which is the case $t=0$) to  freeness (which is the ``limit" when $t$ tends to infinity). 
%Roughly speaking, two subspaces are free if they can be obtained from independent subspaces by conjugating one of them by a Haar unitary element: for all $t\geq 0$,  two subspaces  will be said to be $t$-free if they can be obtained from independent subspaces by conjugating one of them by a  free unitary Brownian motion at time $t$ (with the natural convention that a  free unitary Brownian motion at time $+\infty$ is a Haar unitary element).
We start with the definition of the $t$-free product of two non-commutative \pro spaces. In a few words, it is the space obtained by conjugating one of them, in their tensor product, by a free unitary Brownian motion  at time $t$, free with the tensor product.

Fix $t\in [0,+\infty]$ and let $(\A,\vfi_\A)$ and $(\B,\vfi_\B)$ be  two non-commutative \pro spaces. Let $(\mc{U}^{(t)}, \vfi_{\mc{U}^{(t)}})$ be the   non-commutative \pro space generated by a single unitary element $u_t$ whose distribution is that of a free unitary Brownian motion  at time $t$ (with the   convention that a  free unitary Brownian motion at time $+\infty$ is a Haar unitary element, i.e. a unitary element whose distribution is the uniform law on the unit circle of $\C$).

\begin{definition}[$t$-free product]{The {\em $t$-free product} of $(\A,\vfi_\A)$ and $(\B,\vfi_\B)$, defined up to an isomorphism of non-commutative \pro spaces, is the non-commutative \pro space $(\mc{C},\vfi_{|\mc{C}})$, where $\mc{C}$ is the subalgebra generated by $\A$ and $u_t\B u_t^*$ in %\be\label{20.06.08.2}
$$(\mc{X},\vfi):=[(\A,\vfi_\A)\otimes(\B,\vfi_\B)]*(\mc{U}^{(t)}, \vfi_{\mc{U}^{(t)}}).$$}%\ee and $\vfi_\mc{C}$ is the restriction to $\mc{C}$ of the state of the space of \eqref{20.06.08.2}.
\end{definition}

A few simple observations are in order.

\begin{rmq}\label{20.06.08.14h}{\rm Both $(\A,\vfi_\A)$  and $(\B,\vfi_\B)$ can be identified with subalgebras of the algebra of their $t$-free product (namely with $(\A, \vfi_{|\A})$ and $(u_t\B u_t^*, \vfi_{|u_t\B u_t^*})$).
More specifically, if one defines  
$$\A_{st}:=\{a\in \A\ste \vfi_\A(a)=0, \vfi_\A(aa^*)=1\}, \; \B_{st}:=\{b\in \B\ste \vfi_\B(b)=0, \vfi_\B(bb^*)=1\},$$
then any element in the algebra of the $t$-free product  $(\A,\vfi_\A)$  and $(\B,\vfi_\B)$ can be uniquely written as a constant term plus a linear combination of  words in the elements of $\A_{st}\cup u_t\B_{st}u_t^*$ where no two consecutive letters both belong to $\A_{st} $ or to $u_t\B_{st}u_t^*$.}\end{rmq}
%elements of one of the following types \begin{itemize} \item[(i)] $a_0u_tb_1u_t^*a_1\cdots u_tb_nu_t^*a_n$, with $n\geq 0$, $a_0,\ldots, a_n\in \A_{st}$, $b_1,\ldots, b_b\in \B_{st}$,\item[(ii)] $u_tb_0u_t^*a_1u_tb_1u_t^*\cdots a_n u_tb_nu_t^*$, with $n\geq 0$, $a_1,\ldots, a_n\in \A_{st}$, $b_0,\ldots, b_b\in \B_{st}$, \item[(iii)] $a_1u_tb_1u_t^*\cdots a_n u_tb_nu_t^*$, with $n\geq 1$, $a_1,\ldots, a_n\in \A_{st}$, $b_1,\ldots, b_b\in \B_{st}$, \item[(iv)] $u_tb_1u_t^*a_1\cdots  u_tb_nu_t^*a_n$, with $n\geq 1$, $a_1,\ldots, a_n\in \A_{st}$, $b_1,\ldots, b_b\in \B_{st}$. \end{itemize}
\begin{rmq}{\rm As a consequence, since $u_t$ is unitary and  $(u_t, u_t^*)$ has the same non-commutative distribution as $(u_t^*, u_t)$, the $t$-free product of $(\A,\vfi_\A)$ and $(\B,\vfi_\B)$ is clearly isomorphic, as a non-commutative \pro space, to the $t$-free product of  $(\B,\vfi_\B)$ and $(\A,\vfi_\A)$.}
\end{rmq}

\begin{rmq}\label{20.06.08.16h}{\rm Another consequence of Remark \ref{20.06.08.14h} is that as a  unital algebra,  the algebra of the $t$-free product of $(\A,\vfi_\A)$  and $(\B,\vfi_\B)$ is isomorphic to the free product of the unital algebras $\A/\tilde{\A}$ and $\B/\tilde{\B}$, where $\tilde{\A}$ (resp. $\tilde{\B}$) is the bilateral ideal of the elements $x$ of $\A$ (resp. of $\B$) \st $\vfi_\A(xx^*)=0$ (resp. $\vfi_\B(xx^*)=0$). Thus if $\A$ and $\B$ are subalgebras of the algebra of a non-commutative \pro space $(\mc{M}, \vfi)$, there is a canonical  algebra morphism from  the algebra of the $t$-free product of $(\A,\vfi_{|\A})$  and $(\B,\vfi_{|\B})$ to $\mc{M}$ whose restriction to $\A\cup \B$ preserves the expectation.}
\end{rmq}

Now, we can give the definition of $t$-freeness.  A real $t\in [0,+\infty]$ is still fixed.

\begin{definition}[$t$-freeness] Let $(\mc{M}, \tau)$ be  a non-commutative \pro space. 
\begin{itemize}
\item Two subalgebras $\A, \B$ of $\mc{M}$ are said to be {\em $t$-free} if the canonical algebra morphism   from  the algebra of the $t$-free product of $(\A,\vfi_\A)$  and $(\B,\vfi_\B)$ to $\mc{M}$ mentioned in Remark \ref{20.06.08.16h} preserves the expectation.
\item Two subsets $X,Y$ of $\mc{M}$ are said to be {\em $t$-free} if the subalgebras they generate are $t$-free.
\end{itemize}
\end{definition}

\begin{rmq}{\rm Note that for $t=0$, $t$-freeness is simply the independence, whereas it follows from \cite{haag2} that in the case where $t=+\infty$, it is the freeness.}\end{rmq}

The following proposition is obvious from the definition of $t$-freeness.

\begin{proposition}\label{20.06.08.16h45}
Let $(\mc{M}, \tau)$ be  a non-commutative \pro space. Let $\{a_1,\ldots, a_n\}$ and $\{b_1, \ldots, b_m\}$ be two $t$-free subsets of $\mc{M}$. Then the joint non-commutative distribution of the family $(a_1,\ldots, a_n,b_1, \ldots, b_m)$ depends only on $t$ and on the distributions of the families  $(a_1,\ldots, a_n)$ and  $(b_1, \ldots, b_m)$.
\end{proposition}

\begin{propdef}[Additive and multiplicative $t$-free convolutions] Let us fix $t\in [0,+\infty)$. Let $\mu,\nu$ be compactly supported \pro measures on the real line (resp. on $[0,+\infty)$, on the unit circle). 
Let $a$, $b$ are $t$-free self-adjoint elements  (resp. positive elements, unitary elements) with  distributions $\mu$, $\nu$. Then the distribution of $a+b$ (resp. of $\sqrt{b}a\sqrt{b}$, of $ab$) is a compactly supported \pro measure on the real line which depends only on $t$, $\mu$ and $\nu$, and which will be denoted by $\mu*_t\nu$ (resp. $\mu\odot_t\nu$).\end{propdef}

\begin{pr} Let us treat the case of the sum of two self-adjoint elements. The other cases can be treated analogously. From Proposition \ref{20.06.08.16h45}, it follows that the moments of $a+b$  depend  only on $\mu$ and $\nu$.  To see that these are the moments of a compactly supported \pro measure on the real line,  introduce $M>0$ \st the supports of $\mu$ and $\nu$ are both contained in $[-M,M]$. Then for all $n\geq 1$, by H\"older inequalities in a non-commutative \pro space \cite{Nelson}, $\vfi((a+b)^{2n})\leq 2^{2n}M^{2n}$. By  the remark made after Definition \ref{def ncd}, the result follows.\end{pr}

\begin{proposition}[Matricial model for the $t$-freeness] For each $n\geq 1$, let $M_n$ and $N_n$ be diagonal random matrices whose non-commutative distributions have limits.   Let also, for each $n$, $S_n$ be the matrix of a uniform random permutation of $\{1,\ldots, n\}$ and $U_{n,t}$ be a random $n\times n$ unitary matrix distributed  according to the law of a Brownian motion on the unitary group  at time $t$. Suppose that for each $n$, the sets of random variables $\{M_n,N_n\}, \{S_n\}, \{U_{n,t}\}$ are independent. Then as $n$ tends to infinity, the non-commutative distribution of  $$(M_n, \;U_{n,t}S_nN_nS_n^*U_{n,t}^*)$$ converges in \pro to that of a pair $(a,b)$ of self-adjoint elements of a non-commutative \pro space which are $t$-free.\end{proposition}

\begin{pr}By Theorem \ref{24.06.08.1}, the non-commutative distribution of $(M_n, S_nN_nS_n^*)$ converges to the one of a pair $(x,y)$ of independent elements. Moreover, since for all $n$, the law of $U_n$ is invariant by conjugation, 
by Theorems \ref{24.06.08.1} and \ref{23.06.08.1}, the family of sets $$(\{M_n, S_nN_nS_n^*\}, \{U_{n,t}\})$$ is asymptotically free and the limit distribution of $U_{n,t}$ is that of a  free unitary Brownian motion  at time $t$. By definition of $t$-freeness, this concludes the proof.  \end{pr}

%\begin{corollary}With the notation of the previous proposition, in the case where $M_n,N_n$ have real coefficients and where their limit spectral laws, denoted respectively by $\mu,\nu$ are compactly supported, the measure $\mu*_t\nu$ is the limit spectral law of $M_n+U_{n,t}S_nN_nS_n^*U_{n,t}^*$.
%\end{corollary}

In the next result, the convergences in probability of random measures towards non-random limits are understood with respect to the weak topology on the space of probability measures on the real line.

\begin{corollary}For each $n$,  let $M_n,N_n$ be random $n\times n$ diagonal matrices, one of them having a distribution which is invariant under the action of the symmetric group by conjugation. Suppose that the spectral law of $M_n$ (reps. $N_n$) converges in probability to some compactly supported probability measure $\mu$ (resp. $\nu$) on the real line. Then the spectral law of $M_n+U_{n,t}N_nU_{n,t}^*$ converges in probability to the measure $\mu*_t\nu$.
\end{corollary}

\section{Computation rules for $t$-freeness}\label{10.9.08.2}

\subsection{Multivariate free It\^o calculus}
\subsubsection{Technical preliminaries}In this section, we shall extend some results of \cite{bs98} to the multivariate case. Let us first recall basics of free stochastic calculus. For more involved definitions, the reader should refer to sections 1-2 of \cite{bs98}. Let $(\mc{M},\tau)$ be a faithful\footnote{A non-commutative \pro space $(\mc{M},\tau)$ is said to be {\it faithful} if for all $x$ in $\mc{M}\setminus\{0\}$, $\tau(xx^*)>0$. Any non-commutative \pro space can be quotiented by a bilateral ideal into a faithful space.} non-commutative \pro space endowed with a filtration $(\mc{M}_t)_{t\geq 0}$ and an  $(\mc{M}_t)_{t\geq 0}$-free additive Brownian motion $(X_t)_{t\geq 0}$.  Let $\mc{M}^{op}$ be the opposite algebra of $\mc{M}$ (it is the same vector space, but it is endowed with the product $a\times_{op} b=ba$). We shall denote by $\sharp$ the left actions of the algebra $\mc{M}\otimes \mc{M}^{op}$ on $\mc{M}$ and $\mc{M}\otimes \mc{M}$ defined by $(a\otimes b)\sharp u=aub$ and $(a\otimes b)\sharp (u\otimes v)=au\otimes vb$. The algebras $\mc{M}$ and $\mc{M}\otimes \mc{M}^{op}$ are endowed with the inner products defined by $\lan a,b\ran=\tau(ab^*)$ and $\lan a\otimes b, c\otimes d\ran = \tau(ac^*)\tau(bd^*)$. The Riemann integral of  functions defined on a closed interval with left and right limits at any point with values in the Hilbert space\footnote{The Hilbert space $L^2(\mc{M}, \tau)$ is the completion of $\mc{M}$ with respect to the inner product $\lan x,y\ran=\tau(xy^*)$.}   $L^2(\mc{M}, \tau)$ 
 is  a well known notion.   Now, we shall recall the definition of  the stochastic integral.  
A simple adapted  biprocess is a piecewise constant map $U$ from $[0,+\infty)$ to   $\mc{M}\otimes \mc{M}^{op}$ vanishing for $t$ large enough \st $U_t\in \mc{M}_t\otimes \mc{M}_t$ for all $t$. The set of simple biprocesses   is endowed with the inner product 
$$\lan U,V\ran=\int_0^\infty \lan U_t,V_t\ran\ud t.$$
We shall denote by $\B_2^a$  the closure of the set of simple adapted biprocesses with respect to this inner product.
Let $U$   be a simple adapted biprocess. Then there exists times $0=t_0\leq t_1\leq \cdots \leq t_m$ \st $L$ (resp. $U$) is constant on each $[t_{i}, t_{i+1})$ and vanishes on $[t_m, +\infty)$. Then we define $$ \int_0^\infty U_t\ud X_t=\sum_{i=0}^{m-1}U_{t_i}\sharp(X_{t_{i+1}}-X_{t_i}).$$
It can be proved (Corollary 3.1.2 of \cite{bs98}) that the map $U\mapsto \int_0^\infty U_t\ud X_t$ can be extended isometrically from $\B_2^a$ to  $L^2(\mc{M}, \tau)$.  
   
\subsubsection{Free It\^o processes}
We shall call a {\it free It\^o process} any process 
\be\label{9.9.08.1}A_t=A_0+\int_0^tL_s\ud s+\int_0^t U_s\ud X_s,\ee 
where $A_0\in \mc{M}_0$, $L$ is an adapted   process with left and right limit at any point and $U\in \B_2^a$. In this case, we shall denote \be\label{9.9.08.2}\ud A_t=L_t\ud t+U_t\sharp\ud X_t.\ee 
The part $U_t\sharp \ud X_t$ of this expression is called the {\it martingale part} of $A$. Note that the process $A$ is determined by $A_0$ and $\ud A_t$.  % i.e. for any adapted processes $(C_t), (D_t)$, for all $t$,  \be\label{9.9.08.2}\int_0^t C_s(\ud A_s)D_s=\int_0^tC_sL_sD_s\ud s+\int_0^t (C_s\otimes D_s)\sharp U_s\sharp \ud X_s.\ee

We shall use the following lemma, which follows from Proposition 2.2.2 of \cite{bs98} and from the linearity of $\tau$.
\begin{lemma}\label{martpart=zero}Let $A_t$ be as in \eqref{9.9.08.1}. Then $\tau(A_t)=\tau(A_0)+\int_0^t\tau(L_s)\ud s$.
\end{lemma}

\subsubsection{Multivariate free It\^o calculus}Consider $n$ elements $a_1,\ldots, a_n\in \mc{M}$ for some $n\geq 2$. Consider also two elements %$a_1,\ldots, a_{i-1},a_{i+1}, \ldots, a_{j-1}, a_{j+1}, \ldots, a_n\in \mc{M}$ 
$u=\sum_k x_k\otimes y_k, v=\sum_lz_l\otimes t_l$ of $\mc{M}\otimes \mc{M}^{op}$. For all $1\leq i<j\leq n$, we define an element of $\M$ by setting
\begin{align*}
\lan\lan a_1,\ldots, a_{i-1},u,a_{i+1}, \ldots, a_{j-1},v, a_{j+1}, \ldots, a_n\ran\ran_{i,j}=&\\
& \hskip -3cm  \sum_{k,l}a_1\cdots a_{i-1}x_k\tau(y_ka_{i+1}\cdots a_{j-1}z_l)t_l a_{j+1} \cdots a_n.
\end{align*}

The following theorem follows from Theorem 4.1.12 and the remark following in \cite{bs98}.
\begin{theorem}\label{9.9.08.3}Let $A_t= A_0+\int_0^tL_s\ud s+\int_0^t U_s\ud X_s$ and $B_t= B_0+\int_0^tK_s\ud s+\int_0^t V_s\ud X_s$ be two It\^o processes with respect to the same free Brownian motion $(X_t)$. Then $AB$ is a free It\^o process and with the notations of \eqref{9.9.08.2}, $$\ud(AB)_t=A_t\ud B_t+(\ud A_t)B_t+\lan\lan U_t,V_t\ran\ran_{1,2} \ud t.$$
\end{theorem}

In order to prove computation rules for $t$-freeness, we shall need the following theorem.

\begin{theorem}\label{Ito.multidim}Let $A_1,\ldots, A_n$ be free It\^o processes with respect to the same Brownian motion. For all $k$, denote $A_{k,t}=A_{k,0}+\int_0^tL_{k,s}\ud s+\int_0^t U_{k,s}\ud X_s$. Then $A_1\cdots A_n$    is a free It\^o process and  
\begin{align*}
\ud (A_1\cdots A_n)_t=&\sum_{k=1}^n A_{1,t}\cdots A_{k-1,t}(\ud A_{k,t})A_{k+1,t}\cdots A_{n,t}\\
& +\sum_{1\leq k<l\leq n}\lan\lan   A_{1,t},\ldots, A_{k-1,t},U_{k,t},A_{k+1,t},\ldots, A_{l-1,t},U_{l,t}, A_{l+1,t},\ldots, A_{n,t}\ran\ran_{k,l} \ud t.
\end{align*}
\end{theorem}

\begin{pr}Let us prove this theorem by induction on $n$. For $n=1$, it is obvious. Let us suppose the result to hold at rank $n$. Then the martingale part of $A_1\cdots A_n$ is $$\sum_{k=1}^n A_{1,t}\cdots A_{k-1,t}(U_{k,t}\sharp \ud X_{k,t})A_{k+1,t}\cdots A_{n,t}.$$Thus by Theorem  \ref{9.9.08.3}, $A_1\cdots A_{n+1}$ is a free It\^o process and 
\begin{align*}
\ud(A_1\cdots A_{n+1})_t=& (A_1\cdots A_{n})_t\ud A_{n+1,t}+(\ud (A_1\cdots A_n)_t)A_{n+1,t}\\ 
& +\sum_{k=1}^n \lan\lan A_{1,t},\ldots, A_{k-1,t},U_{k,t},A_{k+1,t},\ldots, A_{n,t},U_{n+1,t}\ran\ran_{k,n+1} \ud t\\
=& \sum_{k=1}^{n+1} A_{1,t}\cdots A_{k-1,t}(\ud A_{k,t})A_{k+1,t}\cdots A_{n,t}\\
&\hskip -2cm +\sum_{1\leq k<l\leq n}\lan\lan   A_{1,t},\ldots, A_{k-1,t},U_{k,t},A_{k+1,t},\ldots, A_{l-1,t},U_{l,t}, A_{l+1,t},\ldots, A_{n,t}, A_{n+1,t}\ran\ran_{k,l} \ud t\\
&\hskip -2cm + \sum_{k=1}^n \lan\lan A_{1,t},\ldots, A_{k-1,t},U_{k,t},A_{k+1,t},\ldots, A_{n,t},U_{n+1,t}\ran\ran_{k,n+1} \ud t,\\
\end{align*}
which concludes the proof.
\end{pr}

\subsection{Computation rules for $t$-freeness}\subsubsection{Main result}
In order to do computations with elements which are $t$-free, we have to find out a formula for the expectation of a product of elements of the type \be\label{gouffran.11.9.08.15h38}x_1u_ty_1u_t^*x_2u_ty_2u_t^*\cdots x_{n}u_ry_nu_t^*,\ee for $\{x_1,\ldots,x_n\}$ independent with $\{y_1,\ldots, y_n\}$ and $\{x_1,y_1,\ldots, x_n,y_n\}$ free   with $u_t$, free unitary Brownian motion. Actually, for the result which follows, the independence of the $x_i $'s and the $y_i$'s will not be useful, thus we consider a non-commutative \pro space $(\mc{M},\tau)$, an integer $n\geq 1$, $a_1,\ldots, a_{2n}\in \M$   and a free unitary Brownian motion $(u_t)$ which is free with $\{a_1,\ldots, a_{2n}\}$. In order to have some more concise formulae, it will be useful to multiply the product of \eqref{gouffran.11.9.08.15h38} by $e^{nt}$. So we define $$f_{2n}(a_1,\ldots,a_{2n},t)=e^{nt}\tau(a_1u_ta_2u_t^*\cdots a_{2n-1}u_ta_{2n}u_t^*).$$
We shall use the convention $f_0(a,t)=\tau(a)$ for all $a\in\M$. 

Since $f_{2n}(a_1,\ldots,a_{2n} ,0)=\tau(a_1\cdots a_{2n})$, the following theorem allows us to deduce all functions $f_{2n}(a_1,\ldots,a_{2n},t)$ (thus the expectation of any product of the type of  \eqref{gouffran.11.9.08.15h38}) from the joint distribution of the $a_i$'s.
%The family of functions $t\mapsto f_{2n}(a_1,\ldots,a_{2n}, t)$ satisfies the following differential system: 

\begin{theorem}\label{suazo.11.9.08} For all $n\geq1$ and all $a_{1}, \ldots, a_{2n}\in \M$ free with the process $(u_t)$, the following differential relations are satisfied:
\begin{align*}
\f{\partial}{\partial t}f_{2n}(a_1,\ldots,a_{2n},t)=&-\!\!\!\!\!\!\!\!\sum_{\substack{1\leq k<l\leq 2n\\ k=l  \; {\rm mod}\;  2}} \!\!\!\!\!\!\!\! 
f_{2n-(l-k)}(a_1,\ldots, a_k,a_{l+1},\ldots, a_{2n},t)f_{l-k}(a_{k+1},\ldots, a_l,t)\\
&\hskip -3cm  +e^t\!\!\!\!\!\!\sum_{\substack{1\leq k<l\leq 2n\\ k\neq l \; {\rm mod}\; 2}} \!\!\!\!\!\!\!\! f_{2n-(l-k)-1}(a_1,\ldots, a_{k-1},a_ka_{l+1},a_{l+2},\ldots, a_{2n} ,t)f_{l-k-1}(a_la_{k+1},a_{k+2},\ldots, a_{l-1} ,t).\\
\end{align*}
\end{theorem}

\begin{pr}Let us introduce the process $(v_t)$ defined by $v_t=e^{t/2}u_t$ for all $t$. As explained in the beginning of section 2.3 of \cite{b97}, this process can be realized as an It\^o process, with the formula $$v_t=1+i\int_0^t \ud X_s v_s.$$ Thus one can realize the family of non-commutative random variables $a_1,\ldots, a_{2n}, (v_t)_{t\geq 0}$ in a faithful non-commutative \pro space  $(\M,\tau)$ endowed with a filtration $(\M_t)_{t\geq 0}$ and an additive free Brownian motion $(X_t)_{t\geq 0}$ \st $a_1,\ldots,a_{2n}\in \M_0$ and for all $t$, $v_t=1+i\int_0^t \ud X_s v_s$ and $v_t^*=1-i\int_0^t v_s^*\ud X_s $.
By definition of $f_{2n}(a_1,\ldots, a_{2n},t)$, one has $$f_{2n}(a_1,\ldots, a_{2n},t)=\tau(a_1v_ta_2v_t^*\cdots a_{2n-1}v_ta_{2n}v_t^*).$$Note that since all $a_i$'s belong to $\M_0$, the processes $A_{1}:=(a_1v_t), A_{2}:=(a_2v_t^*), \ldots,A_{2n-1}:= (a_{2n-1}v_t), A_{2n}:= (a_{2n}v_t^*)$ are all free It\^o processes: if one defines   $U_{k,t}=a_k\otimes iv_t$ for $k$ odd and  $U_{k,t}=-ia_kv_t^*\otimes 1 $ for $k$ even, then  for all $k$, $\ud A_{k,t}=U_{k,t}\sharp \ud X_t$. Thus by theorem \ref{Ito.multidim}, $A_1\cdots A_{2n}$ is an It\^o process \st for all $t$, 
\begin{align*}
(A_1\cdots A_{2n})_t=&(A_1\cdots A_{2n})_0+\int_0^t \sum_{k=1}^{2n} A_{1,s}\cdots A_{k-1,s}(U_{k,s}\sharp \ud X_s)A_{k+1,s}\cdots A_{2n,s} \\
&\hskip -1cm +\int_0^t \sum_{1\leq k<l\leq 2n} \lan\lan A_{1,s},\ldots, A_{k-1,s},U_{k,s},A_{k+1,s},\ldots, A_{l-1,s},U_{l,s}, A_{l+1,s},\ldots, A_{2n,s}\ran\ran_{k,l} \ud s.\\
\end{align*}

Hence by lemma \ref{martpart=zero}, for all $t$,   
\begin{align} \label{cold.war.kids.11.09.08}
\f{\partial}{\partial t}f_{2n}(a_1,\ldots,a_{2n},t)=& \\
&\hskip -2cm  \!\!\!\!\sum_{1\leq k<l\leq 2n}\!\!\!\! \tau(\lan\lan A_{1,t},\ldots, A_{k-1,t},U_{k,t},A_{k+1,t},\ldots, A_{l-1,t},U_{l,t}, A_{l+1,t},\ldots, A_{2n,t}\ran\ran_{k,l}). \nonumber
\end{align}
Now, fix $1\leq k<l\leq 2n$ and discuss according to the parity of $k$ and $l$.

$\bullet$ If $k=l\mod 2$. Suppose for example that $k,l$ are both odd (the other case can be treated in the same way). Then $U_{k,t}=a_k\otimes iv_t$ and $U_{l,t}=a_l\otimes iv_t$, which implies that 
\begin{align*}
\tau(\lan\lan A_{1,t},\ldots, A_{k-1,t},U_{k,t},A_{k+1,t},\ldots, A_{l-1,t},U_{l,t}, A_{l+1,t},\ldots, A_{2n,t}\ran\ran_{k,l})= &\\
&\hskip -8cm i\tau(a_1v_ta_2v_t^*\cdots a_{k-1}v_t^*a_kv_ta_{l+1}v_t\cdots a_{2n}v_t^*)i\tau(v_ta_{k+1}v_t^*\cdots a_{l-1}v_t^*a_l).
\end{align*}
Note that since $\tau$ is tracial and the joint distribution of $a_1,\ldots, a_{2n}, (v_t)_{t\geq 0}$ is the same as the one of $a_1,\ldots, a_{2n}, (v_t^*)_{t\geq 0}$, we have $\tau(v_ta_{k+1}v_t^*\cdots a_{l-1}v_t^*a_l)=\tau(a_{k+1}v_t\cdots a_{l-1}v_ta_lv_t^*)$. Hence 
\begin{align} \label{cocorosie.11.9.08}
\tau(\lan\lan A_{1,t},\ldots, A_{k-1,t},U_{k,t},A_{k+1,t},\ldots, A_{l-1,t},U_{l,t}, A_{l+1,t},\ldots, A_{2n,t}\ran\ran_{k,l})=&\\
&\hskip -8cm -f_{2n-(l-k)}(a_1,\ldots, a_k,a_{l+1},\ldots, a_{2n},t)f_{l-k}(a_{k+1},\ldots, a_l,t).\nonumber
\end{align}

$\bullet$ If $k\neq l\mod 2$. Suppose for example $k$ to be odd and $l$ to be even (the other case can be treated in the same way). Then $U_{k,t}=a_k\otimes iv_t$ and $U_{l,t}=-a_liv_t^*\otimes 1$, which implies that 
\begin{align*}
\tau(\lan\lan A_{1,t},\ldots, A_{k-1,t},U_{k,t},A_{k+1,t},\ldots, A_{l-1,t},U_{l,t}, A_{l+1,t},\ldots, A_{2n,t}\ran\ran_{k,l})=&\\
&\hskip -9cm \tau(a_1v_ta_2v_t^*\cdots a_{k-1}v_t^*a_ka_{l+1}v_t\cdots a_{2n}v_t^*)(-i^2)\tau(v_ta_{k+1}v_t^*\cdots a_{l-1}v_ta_lv_t^*).
\end{align*}
Note that since $v_t^*v_t=e^t$, $\tau$ is tracial and the joint distribution of $a_1,\ldots, a_{2n}, (v_t)_{t\geq 0}$ is the same as the one of $a_1,\ldots, a_{2n}, (v_t^*)_{t\geq 0}$, we have $\tau(v_ta_{k+1}v_t^*\cdots a_{l-1}v_ta_lv_t^*)=e^t\tau(a_la_{k+1}v_t\cdots a_{l-1}v_t^*)$. Hence 
\begin{align}\label{cocorosie.bis.11.9.08}
\tau(\lan\lan A_{1,t},\ldots, A_{k-1,t},U_{k,t},A_{k+1,t},\ldots, A_{l-1,t},U_{l,t}, A_{l+1,t},\ldots, A_{2n,t}\ran\ran_{k,l})=&\\
&\hskip -11cm e^tf_{2n-(l-k)-1}(a_1,\ldots, a_{k-1},a_ka_{l+1},a_{l+2},\ldots, a_{2n} ,t)f_{l-k-1}(a_la_{k+1},a_{k+2},\ldots, a_{l-1} ,t).\nonumber
\end{align}

Equations \eqref{cold.war.kids.11.09.08}, \eqref{cocorosie.11.9.08} and \eqref{cocorosie.bis.11.9.08} together conclude the proof.
\end{pr}\\

The following proposition, which we shall use later, is a consequence of the previous theorem.
\begin{proposition}\label{gulf.shores.13.09.08}
In a non-commutative \pro space $(\mc{M},\tau)$, consider two independent  normal elements $a,b$   with symmetric compactly supported laws. Let $(u_t)$ be a free unitary Brownian motion which is free with $\{a,b\}$. Then the function  $$G(t,z)=\sum_{n\geq 1}\tau((au_tbu_t^*)^{2n})e^{2nt}z^n$$ is the only solution, in a neighborhood of $(0,0)$ in $[0,+\infty) \times \C$,  to the nonlinear, first order partial differential equation \bea \f{\partial G}{\partial t}+4zG  \f{\partial G}{\partial z}&=&0\label{allaround.13.09.08}\\ G(0,z)&=& \sum_{n\geq 1}\tau(a^{2n})\tau(b^{2n})z^n.\label{allaround.13.09.08.1}\eea
\end{proposition}

\begin{pr}Let us define, for all $n\geq 1$, $g_{n}(t)=\tau((au_tbu_t^*)^{n})e^{nt}$. For $n=0$, we set $g_0(t)=0$.  Let us fix $n\geq 1$. In order to apply  the previous theorem, let us define, for $i=1,\ldots, 2n$, $a_i=a$ if $i$ is odd and $a_i=b$ if $i$ is even. By the previous theorem, for all $n\geq 1$, we have 
\begin{align} \label{13.09.08.1}
\f{\partial}{\partial t}g_{n}(t)=&-\!\!\!\!\!\!\!\!\sum_{\substack{1\leq k<l\leq 2n\\ k=l \; {\rm mod}\; 2}} \!\!\!\!\!\!  g_{n-(l-k)/2}(t)g_{(l-k)/2}(t)\\ 
&\hskip -1.5cm +e^t\!\!\!\!\!\!\sum_{\substack{1\leq k<l\leq 2n\\ k\neq l \; {\rm mod}\; 2}} \!\!\!\!\!\! f_{2n-(l-k)-1}(a_1,\ldots, a_{k-1},a_ka_{l+1},a_{l+2},\ldots, a_{2n} ,t)f_{l-k-1}(a_la_{k+1},a_{k+2},\ldots, a_{l-1} ,t).\nonumber
\end{align}

Now, note that since for any $\eps, \eps'=\pm 1$, the joint distribution of $(a,b, u_t)$ is the same as the one of $(\eps a,\eps' b,u_t)$, $g_p(t)=0$ when $p$ is odd. Thus in the first sum of \eqref{13.09.08.1} only pairs $(k,l)$ \st $k=l\mod 4$ have a non null contribution. For the same reason, all terms in the second sum are null. Indeed,  for any $1\leq k<l\leq 2n$ \st $k\neq l \mod 2$,  the set  $\{k+1, k+2, \ldots, l\}$, whose cardinality is odd, has either an odd number of odd elements or an odd number of even elements. To sum up, for all $n\geq 1$, we have 
\begin{align*}
\frac{\partial}{\partial t}g_{2n}(t) &= - \sum_{\substack{1\leq k<l\leq 4n \\ k=l \; {\rm mod} \; 4}} g_{2n-(l-k)/2}(t)g_{(l-k)/2}(t)\\
&= - 4 \sum_{i=1}^{n-1} (n-i) g_{2(n-i)}(t) g_{2i}(t)\\
&= - 2n \sum_{i=1}^{n-1}  g_{2(n-i)}(t) g_{2i}(t).
\end{align*}

%\begin{align*}
%\f{\partial}{\partial t}g_{2n}(t)&=-\!\!\!\!\!\!\!\!\sum_{\substack{1\leq k<l\leq 4n\\ k=l \mod 4}} \!\!\!\!\!\!\!\! g_{2n-(l-k)/2}(t)g_{(l-k)/2}(t) = -4\!\!\!\!\sum_{\substack{1\leq i<j\leq n}} \!\!\!\! 
%g_{2(n-(j-i))}(t)g_{2(j-i)}(t)\\ 
%&\hskip -1cm=  -4\sum_{m=1}^{n-1}mg_{2m}(t)g_{2(n-m)}(t)= -2\sum_{m=1}^{n-1}mg_{2m}(t)g_{2(n-m)}(t)-2\sum_{m=1}^{n-1}(n-m)g_{2m}(t)g_{2(n-m)}(t)\\
%&\hskip -1cm = -2n\sum_{m=1}^{n-1}g_{2m}(t)g_{2(n-m)}(t) =-2n\sum_{m=0}^{n}g_{2m}(t)g_{2(n-m)}(t).
%\end{align*}

Thus since $g_0(t)=0$ and $G(t,z)=\sum_{n\geq 1}g_{2n}(t)z^n=\sum_{n\geq 0}g_{2n}(t)z^n$, the last computation implies 
$$\f{\partial G}{\partial t}=-2z\f{\partial G^2}{\partial z},$$
which proves  \eqref{allaround.13.09.08}. The formula \eqref{allaround.13.09.08.1} is obvious. 

To prove the uniqueness, let $H(t,z)=\sum_{n\geq 0}h_n(t)z^n$ be another solution of   \eqref{allaround.13.09.08} and \eqref{allaround.13.09.08.1}. By \eqref{allaround.13.09.08.1}, for all $n\geq 0$, we have $h_n(0)=g_{2n}(0)$ and by  \eqref{allaround.13.09.08.1}, for all $n\geq 0$, we have $\f{\partial}{\partial t}h_{n}(t)=
-2n\sum_{m=0}^{n}h_{n-m}(t)h_{m}(t)$, which implies that $h_0=0$ and that by induction on $n$, $h_n=g_{2n}$.
\end{pr}

\subsubsection{Examples} Let us give examples of applications of the computation rules that we have just established. The third example below is a rather big formula, but we shall need it when we study the problem of existence of $t$-free cumulants. So, let $\A$ and $\B$ be two independent subalgebras of a non-commutative \pro space $(\M,\tau)$ and $(u_t)$ be a free unitary Brownian motion free from $\A\cup \B$. 

1) For $a\in \A, b\in \B$, for all $t\geq 0$, we have \be\label{eddie.vedder.11.9.08}\tau(au_tbu_t^*)=\tau(a)\tau(b).\ee(In fact, it even follows from theorem \ref{suazo.11.9.08} that without the assumption that $a$ and $b$ are independent, for all $t$, we have $\tau(au_tbu_t^*)=e^{-t}\tau(ab)+(1-e^{-t})\tau(a)\tau(b)$).

2) For $a,a'\in \A, b,b'\in \B$, for all $t\geq 0$, we have 
\be\label{eddie.vedder.11.9.08.bis}\tau(au_tbu_t^* a'u_tb'u_t^*)=\ee $$(\tau(a)\tau(a')\tau(bb')+\tau(aa')\tau(b)\tau(b')-\tau(a)\tau(a')\tau(b)\tau(b'))(1-e^{-2t})+\tau(aa')\tau(bb')e^{-2t}.$$

3) For $a,a',a''\in \A$ and $b,b',b''\in \B$, we have
\begin{eqnarray}\label{eddie.vedder.11.9.08.third}
&&\hskip 1cm \tau(au_tbu_t^*a'u_tb'u_t^*a''u_tb''u_t^*)\\ &=&\tau(a)\tau(a')\tau(a'')\tau(b)\tau(b')\tau(b'')(2-6e^{-2t}+4e^{-3t})\nonumber\\
&&-(1-3e^{-2t}+2e^{-3t})\tau(a)\tau(a')\tau(a'')\tau(b)\tau(b'b'')\nonumber\\
&&-(1-3e^{-2t}+2e^{-3t})\tau(a)\tau(a')\tau(a'')\tau(bb')\tau(b'')\nonumber\\
&&-(1-3e^{-2t}+2e^{-3t})\tau(a)\tau(a')\tau(a'')\tau(bb'')\tau(b')\nonumber\\
&&-(1-3e^{-2t}+2e^{-3t})
\tau(aa')\tau(a'')\tau(b)\tau(b')\tau(b'')\nonumber\\
&&-(1-3e^{-2t}+2e^{-3t})\tau(aa'')\tau(a')\tau(b)\tau(b')\tau(b'')\nonumber\\
&&-(1-3e^{-2t}+2e^{-3t})\tau(a)\tau(a'a'')\tau(b)\tau(b')\tau(b'')\nonumber\\
&&+(1-3e^{-2t}+2e^{-3t})[\tau(a)\tau(a')\tau(a'')\tau(bb'b'')+\tau(aa'a'')\tau(b)\tau(b')\tau(b'')]\nonumber\\
&&-(e^{-2t}-e^{-3t}) [\tau(aa')\tau(a'')\tau(bb')\tau(b'')+\tau(aa')\tau(a'')\tau(bb'')\tau(b')+\tau(aa'')\tau(a')\tau(bb'')\tau(b')\nonumber\\ &&+\tau(aa'')\tau(a')\tau(b)\tau(b'b'')+\tau(a)\tau(a'a'')\tau(b'')\tau(bb')+\tau(a)\tau(a'a'')\tau(b)\tau(b'b'')]\nonumber\\
&&+(1-2e^{-2t}+e^{-3t}) [\tau(aa')\tau(a'')\tau(b'b'')\tau(b)+\tau(aa'')\tau(a')\tau(bb')\tau(b'')+\tau(a)\tau(a'a'')\tau(bb'')\tau(b')]\nonumber\\
&&+(e^{-2t}-e^{-3t})\tau(bb'b'')[\tau(a)\tau(a'a'')+\tau(aa')\tau(a'')+\tau(aa'')\tau(a')]\nonumber\\
&&+(e^{-2t}-e^{-3t})\tau(aa'a'')[\tau(b)\tau(b'b'')+\tau(bb')\tau(b'')+\tau(bb'')\tau(b')]\nonumber\\ 
&&+e^{-3t}\tau(aa'a'')\tau(bb'b'')\nonumber
\end{eqnarray}

It can be verified that the last formula actually corresponds to the formula of $\E(aba'b'a''b'')$ with $\{a,a',a''\}$ and $\{b,b',b''\}$ independent when $t=0$, and to the formula of $\tau(aba'b'a''b'')$ with  $\{a,a',a''\}$ and $\{b,b',b''\}$ free when $t$ tends to infinity.

\subsection{Multiplicative and additive $t$-free convolutions of two symmetric Bernoulli laws}
\subsubsection{Multiplicative case} 
Here, we shall  compute the multiplicative $t$-free convolution of $\f{\delta_{-1}+\delta_1}{2}$ (considered as a law on the unit circle) with itself.
\begin{theorem}For all $t\geq 0$, $\f{\delta_{-1}+\delta_1}{2}\odot_t\f{\delta_{-1}+\delta_1}{2}$ is the only law on the unit circle which is invariant under the symmetries with respect to the real and imaginary axes and whose push-forward by the map $z\mapsto z^2$ has the law of $u_{4t}$, a free unitary Brownian motion taken at time $4t$. 
\end{theorem}

\begin{remark}{\rm The moments of $u_{4t}$ have been computed by P. Biane at Lemma 1 of \cite{b97}: for all $n\geq 1$, 
\begin{equation}\label{moments mbul}
\tau(u_{4t}^n)=\f{e^{-2nt}}{n}\sum_{k=0}^{n-1}\f{(-4nt)^k}{k!}\binom{n}{k+1}.
\end{equation}}
\end{remark}

\begin{pr} In a non-commutative \pro space $(\mc{M},\tau)$, consider two independent  normal   elements $a,b$   with law $\f{\delta_{-1}+\delta_1}{2}$. Let $(u_t)$ be a free unitary Brownian motion which is free with $\{a,b\}$. Then $\f{\delta_{-1}+\delta_1}{2}\odot_t\f{\delta_{-1}+\delta_1}{2}$ is the distribution of the unitary element $au_tbu_t^*$. Since the joint distribution of $(a,b,u_t)$ is the same as the one of $(-a,b ,u_t)$, $\f{\delta_{-1}+\delta_1}{2}\odot_t\f{\delta_{-1}+\delta_1}{2}$ is invariant under the transformation $z\mapsto -z$. Moreover, $(au_tbu_t^*)^*=u_tbu_t^*a$ has the same distribution as $au_tbu_t^*$ (because $\tau$ is tracial and $u_t$ has the same law as $u_t^*$), hence $\f{\delta_{-1}+\delta_1}{2}\odot_t\f{\delta_{-1}+\delta_1}{2}$ is invariant under the transformation $z\mapsto \bar{z}$. This proves that $\f{\delta_{-1}+\delta_1}{2}\odot_t\f{\delta_{-1}+\delta_1}{2}$  is invariant under the symmetries with respect to the real and imaginary axes.  

Since any distribution on the unit circle is determined by its moments, to prove that the push-forward of  $\f{\delta_{-1}+\delta_1}{2}\odot_t\f{\delta_{-1}+\delta_1}{2}$ by   the map $z\mapsto z^2$ is the law of $u_{4t}$, it suffices to prove that for all $n\geq 1$, $$\tau((au_tbu_t^*)^{2n})=\tau(u_{4t}^{n}),$$ i.e. to prove that the functions $$F_1(t,z)=\sum_{n\geq 1} \tau((au_tbu_t^*)^{2n})e^{2nt}z^n\quad\textrm{ and  }\quad F_2(t,z)=\sum_{n\geq 1} \tau(u_{4t}^{n})e^{2nt}z^n$$ are equal. It follows from Proposition \ref{gulf.shores.13.09.08} that $F_1$  is the only solution, in a neighborhood of $(0,0)$ in $[0,+\infty) \times \C$, to equation \eqref{allaround.13.09.08} satisfying $F_1(0,z)=\f{z}{1-z}$. But it follows from Lemma 1 of \cite{b97} that $F_2$ is also a solution of \eqref{allaround.13.09.08} with the same initial conditions. By uniqueness, it closes the proof.\end{pr}

For all $t\in [0,1]$, let us define $\beta(t)=2\sqrt{t(1-t)} +\arccos(1-2t)$. Then $\beta(t)$ is an increasing function of $t$ which goes from $0$ to $\pi$ when $t$ goes from $0$ to $1$. P. Biane has proved in \cite[Prop. 10]{b97b} that the distribution of $u_{4t}$ is absolutely continuous with respect to the Lebesgue measure on the unit circle, that its support is the full unit circle for $t\geq 1$, and the set $\{e^{i\theta} : |\theta|\leq \beta(t)\}$ for all $t\in [0,1]$. Moreover, the density of this distribution with respect to the uniform probability measure on the unit circle, which we denote by $\rho_{4t}$, is positive and analytic on the interior of its support for all $t\geq 0$, except at $-1$ for $t=1$.

\begin{figure}[h!]
\begin{center}
\rotatebox{270}{\scalebox{0.6}{\includegraphics{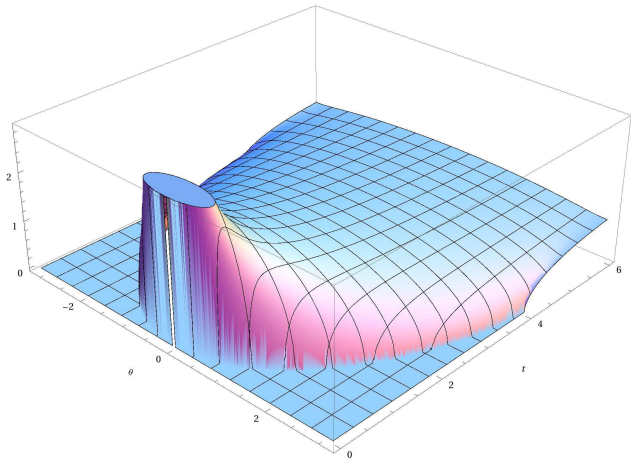}}}
\caption{The density of the distribution of $u_t$ at the point $e^{i\theta}$ in function of $\theta$ and $t$. One sees the support progressively filling the circle when $t$ increases from $0$ to $4$, and then the distribution rapidly converging towards the uniform distribution.}
\end{center}
\end{figure}

\begin{remark}\label{26.11.08.11h26}{\rm  Since there is no simple formula for the density of $u_{t}$, it may be worth explaining how we got the picture above. The expression of the moments of the distribution of $u_{t}$ given by (\ref{moments mbul}) is numerically ineffective, because it is an alternated sum of very large numbers. It only allows one to compute the first few dozens of moments of the distribution. Nevertheless, this expression of the moments allows one to prove that, for all $t\geq 0$, the function $\kappa_{t}$ defined on the interior of the complex unit disk by the formula
$$\kappa_{t}(z)=\int_{0}^{2\pi} \frac{e^{i\theta}+z}{e^{i\theta}-z} \rho_{t}(e^{i\theta})\; \frac{d\theta}{2\pi}=1+2 \sum_{k=1}^{+\infty} \tau(u_{t}^{k})z^{k},$$
satisfies the equation
\begin{equation}\label{calc dens}
\frac{\kappa_{t}(z)-1}{\kappa_{t}(z)+1}e^{\frac{t}{2}\kappa_{t}(z)}=z.
\end{equation}
This fact can be established using the Lagrange inversion formula (see \cite{b97}), see also \cite[4.2.2]{b97b}. Now, on one hand, a computer seems to be able to solve this equation more reliably than it computes the moments of the distribution. On the other hand, $\kappa_{t}$ is the holomorphic function in the unit disk whose real part is the harmonic extension of the density of the distribution of $u_{t}$. Thus, we evaluated $\rho_{t}(e^{i\theta})$ by taking the real part of a numerical solution of (\ref{calc dens}) with $z=e^{i\theta}$.
}
\end{remark}

From the facts exposed above Remark \ref{26.11.08.11h26}, one deduces easily the next result. 

\begin{corollary}\label{support mult} For all $t>0$, the measure $\f{\delta_{-1}+\delta_1}{2}\odot_t\f{\delta_{-1}+\delta_1}{2}$ 
has a density with respect to the uniform probability measure on the unit circle, which we shall denote by $\sigma_t$ and which is given by the formula $\sigma_t(z)=\rho_{4t}(z^2)$ for all $z$ in the unit circle. In particular, the support of this measure is
the full unit circle for $t\geq 1$ and the set $\{e^{i\theta} : |\theta|\leq \frac{1}{2} \beta(t) \mbox{ or }  |\pi- \theta|\leq \frac{1}{2} \beta(t)\}$ for $t\in [0,1]$. The density $\sigma_t$ is positive and analytic on the interior of its support for all $t\geq 0$, except at $\pm i$ for $t=1$.  
\end{corollary}

\begin{figure}[h!]
\begin{center}
\rotatebox{270}{\scalebox{0.6}{\includegraphics{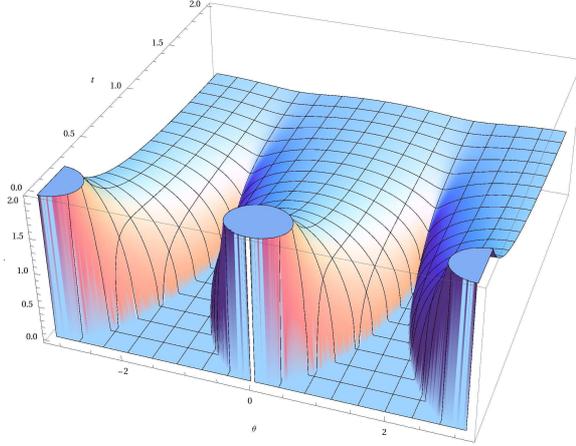}}}
\caption{The density of  $\f{\delta_{-1}+\delta_1}{2}\odot_t\f{\delta_{-1}+\delta_1}{2}$ at the point $e^{i\theta}$ in function of $\theta$ and $t$. The support progressively fills the circle when $t$ increases from $0$ to $1$, and then the distribution converges rapidly towards the uniform distribution.}
\end{center}
\end{figure}

\subsubsection{Additive case} 
Here, we shall  compute the additive $t$-free convolution of $\f{\delta_{-1}+\delta_1}{2}$ (considered as a law on the real line) with itself.
\begin{theorem}For all $t\geq 0$, $\f{\delta_{-1}+\delta_1}{2}*_t\f{\delta_{-1}+\delta_1}{2}$ is the only symmetric law on the real line whose push-forward by the map $x\mapsto x^2$ has the law of $2+v+v^*$, with $v$ unitary element distributed according to $\f{\delta_{-1}+\delta_1}{2}\odot_t\f{\delta_{-1}+\delta_1}{2}$. 
\end{theorem}

%{\bf Thierry: } un truc a faire eventuellement ici est de donner une expression simple pour les moments de $\f{\delta_{-1}+\delta_1}{2}*_t\f{\delta_{-1}+\delta_1}{2}$, voire de sa densite. J'ai essayé, je n'y suis pas arrive.

\begin{pr} In a non-commutative \pro space $(\mc{M},\tau)$, consider two independents  normal   elements $a,b$   with law $\f{\delta_{-1}+\delta_1}{2}$. Let $(u_t)$ be a free unitary Brownian motion which is free with $\{a,b\}$. Then $\f{\delta_{-1}+\delta_1}{2}*_t\f{\delta_{-1}+\delta_1}{2}$ is the distribution of $a+u_tbu_t^*$. Since the joint distribution of $(a,b,u_t)$ is the same as the one of $(-a,-b ,u_t)$, $\f{\delta_{-1}+\delta_1}{2}*_t\f{\delta_{-1}+\delta_1}{2}$ is symmetric.  Note that since $a^2$ and $b^2$ have $\delta_1$ for distribution,  one can suppose that $a^2=b^2=1$. In this case, $$(a+u_tbu_t^*)^2=2+a u_tbu_t^*+ u_tbu_t^*a=2+a u_tbu_t^*+(a u_tbu_t^*)^*,$$and the result is obvious by definition of $\odot_t$.
\end{pr}

From the last result and Proposition \ref{support mult}, we deduce the following.

\begin{corollary} For all $t>0$, the measure $\f{\delta_{-1}+\delta_1}{2}*_t \f{\delta_{-1}+\delta_1}{2}$ 
has a density with respect to the Lebesgue measure on $[-2,2]$, which we shall denote by $\eta_t$ and which is given by the formula
$$\forall x\in [-2,2] \; , \;\; \eta_t(x)= \rho_{4t}(e^{4i\arccos \frac{x}{2}}) \frac{1}{\pi\sqrt{4-x^2}}.$$
The support of this measure is the interval $[-2,2]$ for $t\geq 1$, and the set
$$\left[-2,-2 \cos \frac{\beta(t)}{4}\right] \cup \left[-2 \sin \frac{\beta(t)}{4},2 \sin \frac{\beta(t)}{4}\right] \cup \left[2 \cos \frac{\beta(t)}{4}, 2\right]$$
for $t\in [0,1]$.
The density $\eta_t$ is positive and analytic on the interior of its support for all $t\geq 0$, except at $\pm \sqrt{2}$ for $t=1$. 
\end{corollary}

\begin{figure}[h!]
\begin{center}
\rotatebox{270}{\scalebox{0.6}{\includegraphics{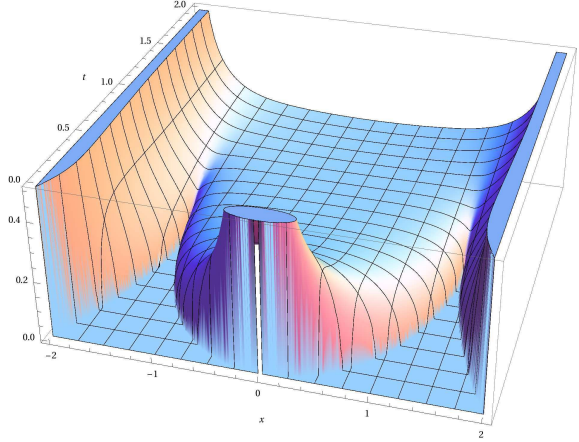}}}
\caption{The density of  $\f{\delta_{-1}+\delta_1}{2}*_t\f{\delta_{-1}+\delta_1}{2}$ at the point $x$ in function of $x$ and $t$. The support fills progressively the interval $[-2,2]$ when $t$ increases from $0$ to $1$, and then the distribution converges rapidly towards the arcsine distribution.}
\end{center}
\end{figure}

%\begin{remark}{\rm One recovers the fact that as $t$ goes from $0$ to $+\infty$, $\f{\delta_{-1}+\delta_1}{2}*_t\f{\delta_{-1}+\delta_1}{2}$  passes from $\f{\delta_{-1}+\delta_1}{2}*\f{\delta_{-1}+\delta_1}{2}=\f{1}{4}(\delta_{-2}+2\delta_0+\delta_2)$ to $\f{\delta_{-1}+\delta_1}{2}\bxp\f{\delta_{-1}+\delta_1}{2}= 1_{[-2,2]}(x) \f{\ud x}{\pi \sqrt{2-x^2}}$.}\end{remark}

\section{The lack of cumulants}\label{10.9.08.3}

In this section, we investigate the existence of an analogue of classical and free cumulants in the context of $t$-freeness. 
Informally, the problem is to find multilinear forms defined on any non-commutative probability space which vanish when evaluated on a family of elements which can be split into two non-empty subfamilies which are $t$-free. 

More precisely, given a non-commutative probability space $(\mc{M},\varphi)$, we would like to know if there exists a family $(k_{n})_{n\geq 2}$ of multilinear forms on $\mc{M}$, with $k_n$ an $n$-linear form for all $n\geq 2$, such that, for all $n\geq 2$, all $n_{1},n_{2}>0$ such that $n_{1}+n_{2}=n$, all $m_{1},\ldots,m_{n}$ in $\mc{M}$ such that $\{m_{1},\ldots,m_{n_{1}}\}$ and $\{m_{n_{1}+1},\ldots,m_{n_{1}+n_{2}}\}$ are $t$-free, and finally for all $\sigma \in \Sym_{n}$, one has $k_{n}(m_{\sigma(1)},\ldots,m_{\sigma(n)})=0$. 

Our main result is negative: there does not exist in general such a family of multilinear forms, at least in a large class which we describe now.

\begin{definition} Let $(\mc{M},\varphi)$ be a non-commutative probability space. Let $n\geq 1$ be an integer. Let $\sigma$ be an element of $\Sym_{n}$. We define the $n$-linear form $\varphi_{\sigma}$ on $\mc{M}$ as follows:
$$\forall m_{1},\ldots,m_{n} \in \mc{M} \; , \;\; \varphi_{\sigma}(m_{1},\ldots,m_{n})=\prod_{\substack{c {\rm \; cycle \; of\; }\sigma\\c=(i_{1}\ldots i_{r})}} \varphi(m_{i_{1}}\ldots m_{i_{r}}).$$
\end{definition}

Using only the algebra structure of $\mc{M}$ and the linear form $\varphi$, a linear combination of the forms $\{\varphi_\sigma:\sigma\in \Sym_n\}$ seems to be the most general $n$-linear form that one can construct on $\mc{M}$.
We seek cumulants within this wide class of $n$-linear forms. Our definition does not require that the vanishing of cumulants characterize $t$-freeness. We only insist that mixed cumulants of $t$-free variable vanish.

\begin{definition} Let $n\geq 2$ be an integer. Let $t\geq 0$ be a real number. A {\em $t$-free cumulant of order $n$} is a collection $(c(\sigma))_{\sigma\in \Sym_{n}}$ of complex numbers such that $\displaystyle \sum_{\sigma \; n{\rm -cycle}} c(\sigma)\neq 0$ and the following property holds for every non-commutative probability space $(\mc{M},\varphi)$ : for any pair $(\A,\mc{B})$ of 
 sub-algebras of $\mc{M}$ which are $t$-free with respect to $\varphi$, for any family $(m_{1},\ldots,m_{n})$ of elements of $\mc{A}\cup \mc{B}$, which do not all belong to $\mc{A}$, and not all to $\mc{B}$, we have
\begin{equation}\label{cum 1}
\sum_{\sigma\in \Sym_{n}} c(\sigma) \varphi_{\sigma}(m_{1},\ldots,m_{n})=0.
\end{equation}
\end{definition}

Let us emphasize that what we call {\em cumulant} is not a specific multilinear form, but rather a collection of coefficients which allows one to define a multilinear form on any non-commutative probability space.

If $(c(\sigma))_{\sigma\in \Sym_{n}}$ is a $t$-free cumulant of order $n$ and $m_{1},\ldots,m_{n}$ are elements of a non-commutative probability space $(\mc{M},\varphi)$, at least one of which is equal to $1$, then
\begin{equation}\label{cum 2}
\sum_{\sigma\in \Sym_{n}} c(\sigma) \varphi_{\sigma}(m_{1},\ldots,m_{n})=0.
\end{equation}
Indeed, the subalgebra $\C.1$ of $\mc{M}$ is $t$-free with any subalgebra of $\mc{M}$.

We extend the previous definition by including the free case $t=+\infty$. We will mainly consider collections $(c(\sigma))_{\sigma\in \Sym_{n}}$ with the property that $c(\rho\sigma\rho^{-1})=c(\sigma)$ for all $\sigma,\rho \in \Sym_{n}$. We call such collections {\em conjugation-invariant}. They are in fact indexed by conjugacy classes of $\Sym_{n}$, that is, integer partitions of $n$. Thus, we write use as well the notation $(c_{\lambda})_{\lambda \vdash n}$ for a conjugation-invariant collection.

Our main results are the following.

\begin{theorem}\label{existence} For all $t\in [0,+\infty]$ and all $n\in \{2,3,4,5,6\}$, there exists, up to scaling, a unique conjugation-invariant $t$-free cumulant of order $n$.
\end{theorem}

\begin{theorem} \label{non existence} There exists a $t$-free cumulant of order $7$ if and only if $t=0$ or $t=+\infty$.
\end{theorem}

Let us start by proving that we lose nothing by focusing on conjugation-invariant $t$-free cumulants.

\begin{lemma} \label{cum ci} If for some $t$ and some $n$ there exists a $t$-free cumulant of order $n$, then there exists such a cumulant $(c(\sigma))_{\sigma\in \Sym_{n}}$ such that moreover $c(\sigma)=c(\rho \sigma \rho^{-1})$ for all $\sigma,\rho\in \Sym_{n}$.
\end{lemma}

\begin{proof} The point is that the order of the arguments is arbitrary in (\ref{cum 1}). Hence, if (\ref{cum 1}) is satisfied, then for all $\rho\in \Sym_{n}$,
\begin{align*}
0=\sum_{\sigma\in \Sym_{n}} c(\sigma) \varphi_{\sigma}(m_{\rho^{-1}(1)},\ldots,m_{\rho^{-1}(n)})&=\sum_{\sigma\in \Sym_{n}} c(\sigma)\varphi_{\rho^{-1}\sigma\rho}(m_{1},\ldots,m_{n})\\
&=\sum_{\sigma\in \Sym_{n}} c(\rho\sigma\rho^{-1})\varphi_{\sigma}(m_{1},\ldots,m_{n}).
\end{align*}
Hence, if $(c(\sigma))_{\sigma\in \Sym_{n}}$ is a $t$-free cumulant, then so is $(c(\rho\sigma \rho^{-1}))_{\sigma\in \Sym_{n}}$. By averaging over $\rho$, we get a conjugation-invariant cumulant. \end{proof}

Observe that the assumption made in the definition of a cumulant that the sum of $c(\sigma)$ when $\sigma$ spans the $n$-cycles is nonzero implies that $c_n\neq 0$ for any conjugation-invariant cumulant.

%If $\sigma$ is a permuation, we write $[\sigma]=\lambda$ to indicate that the conjugacy class of $\sigma$ is represented by the partition $\lambda$. 
%Let us introduce another notation. On a non-commutative probability space $(\mc{M},\varphi)$, given $\lambda \vdash n$, we set
%$$\varpi_{\lambda}=\sum_{[\sigma]=\lambda} \varphi_{\sigma}.$$
%Thus, a conjugation-invariant cumulant has the form $\sum_{\lambda\vdash n}c_{\lambda} \varpi_{\lambda}$.

Let us introduce some notation. Given a permutation $\sigma$ of $\{1,\ldots,n\}$, we denote by $\{\{\sigma\}\}$ the partition of $\{1,\ldots,n\}$ whose blocks are the sets underlying the cycles of $\sigma$. Let $\mc{P}(n)$ denote the set of partitions of the set $\{1,\ldots,n\}$. Let $(\mc{A},\varphi)$ be a commutative non-commutative probability space. For each partition $\pi\in \mc{P}(n)$, we define an $n$-linear form $\varphi_{\pi}$ on $\mc{A}$ by setting $\varphi_{\pi}=\varphi_{\sigma}$, where $\sigma$ is any permutation of $\{1,\ldots,n\}$ such that $\{\{\sigma\}\}=\pi$. Since $\mc{A}$ is commutative, this definition does not depend on the choice of $\sigma$. Finally, when $\sigma$ is a permutation of $\{1,\ldots,n\}$, we say that $i,j\in \{1,\ldots,n\}$ are {\em consecutive} in a cycle of $\sigma$ if $\sigma(i)=j$ or $\sigma(j)=i$. We will use repeatedly the following fact, which is a consequence of Proposition \ref{20.06.08.16h45} and Proposition 1 of \cite{b98}.

\begin{lemma}\label{c sigma pi pi} Choose two integers $k,l>0$ and set $n=k+l$. \\
1. There exists universal coefficients $(C(\sigma,\pi,\pi'))_{\sigma\in \Sym_n,\pi\in \mc{P}(k),\pi'\in \mc{P}(l)}$ such that the following property holds:

Let $\mc{A}$ and $\mc{B}$ be two commutative sub-algebras of some non-commutative probability space $(\mc{M},\varphi)$ which are $t$-free with respect to $\varphi$.  Consider $\sigma\in \Sym_{n}$. For all $a_{1},\ldots,a_{k}\in \mc{A}$ and all $b_{1},\ldots,b_{l}\in \mc{B}$, 
\begin{equation}\label{expand}
\varphi_{\sigma}(a_{1},\ldots,a_{k},b_{1},\ldots,b_{l})=\sum_{\pi\in \mc{P}(k),\pi'\in \mc{P}(l)} C(\sigma,\pi,\pi') \varphi_{\pi}(a_{1},\ldots,a_{k})\varphi_{\pi'}(b_{1},\ldots,b_{l}).
\end{equation}
2. The coefficient $C(\sigma,\pi,\pi')$ can be non-zero only if every block of $\pi$ is contained in a block of $\{\{\sigma\}\}$.\\
3. If two elements $i$ and $j$ of $\{1,\ldots,k\}$ are consecutive in a cycle of $\sigma$, then $C(\sigma,\pi,\pi')$ can be non-zero only if $i$ and $j$ are in the same block of $\pi$.
\\
4. The parts 2. and 3. of this lemma are also valid for $\pi'$ (modulo a translation of $k$ of the indices, since $\pi'$ is a partition of $\{1,\ldots, l\}$ and not of $\{k+1,\dots, k+l\}$).
\end{lemma}

With the notation of the lemma above, we associate to every collection $(c(\sigma))_{\sigma\in \Sym_{n}}$ the following family of coefficients:
\begin{equation}\label{def D}
\forall \pi\in \mc{P}(k), \pi'\in \mc{P}(l), \;\; D_{c}(\pi,\pi')=\sum_{\sigma\in \Sym_{n}} c(\sigma)C(\sigma,\pi,\pi'),
\end{equation}
which will play an important role in the proofs of Theorems \ref{non existence} and \ref{existence}. \\

\begin{proof} (Theorem \ref{non existence}) Let us choose $t>0$ a positive real. We prove by contradiction that there exists no $t$-free cumulant of order $7$. So, let us assume that there exists one and let $(c(\sigma))_{\sigma\in \Sym_7}$ be one of them, which we choose to be conjugation-invariant thanks to Lemma \ref{cum ci}. Thus, we denote it also by $(c_{\lambda})_{\lambda\vdash 7}$. Since $c_7\neq 0$, we may and will assume that $c_{7}=1$. Then, we proceed as follows.

Let us consider a non-commutative probability space $(\mc{M},\varphi)$ and two commutative sub-algebras $\mc{A}$ and $\mc{B}$ of $\mc{M}$ which are $t$-free with respect to $\varphi$. Let us choose $a_{1},a_{2},a_3\in \mc{A}$ and $b_{1},b_2,b_3,b_4,b_5 \in \mc{B}$, which we assume to be all centered. Set $k_{7}=\sum_{\sigma\in \Sym_{7}} c(\sigma)\varphi_{\sigma}$. By using the $t$-freeness of $\A$ and $\mc{B}$, we will express $k_7(a_1,a_2,b_1,b_2,b_3,b_4,b_5)$ and $k_7(a_1,a_2,a_3,b_1,b_2,b_3,b_4)$ in terms of the coefficients $(c_\lambda)_{\lambda\vdash 7}$, the joint moments of $a_1,a_2,a_3$, and the joint moments of $b_1,b_2,b_3,b_4,b_5$. By the assumption that $k_7$ is a $t$-free cumulant, the two expressions that we get must vanish. Since the joint  distributions of the $a$'s and of the $b$'s are both arbitrary among those of families of centered elements, every coefficient of a given product of moments of the $a$'s and $b$'s must vanish. This gives us linear relations on the coefficients $(c_\lambda)_{\lambda\vdash 7}$, which will turn out to be incompatible.\\

Let us start with $k_{7}(a_{1},a_{2},b_{1},b_2,b_3,b_4,b_{5})$. By Lemma \ref{c sigma pi pi}, this quantity can be written as
\begin{align}\label{dev k7 52}
&\sum_{\sigma\in \Sym_7,\pi\in\mc{P}(2),\pi'\in\mc{P}(5)} c(\sigma) C(\sigma,\pi,\pi') \varphi_\pi(a_1,a_2) \varphi_{\pi'}(b_1,b_2,b_3,b_4,b_5) \nonumber \\
&\hskip 4cm =\sum_{\pi\in\mc{P}(2),\pi'\in\mc{P}(5)} D_{c}(\pi,\pi')\varphi_\pi(a_1,a_2) \varphi_{\pi'}(b_1,b_2,b_3,b_4,b_5).
\end{align}
We are thus interested in computing, for each pair $(\pi,\pi')$, the coefficient $D_{c}(\pi,\pi')$. It turns out to be convenient to think of $b_1,\ldots,b_5$ as occupying the slots $3$ to $7$ rather than $1$ to $5$ and to see $\pi'$ as a partition of the set $\{3,\ldots,7\}$ accordingly. We hope that no ambiguity will result from this change in our conventions.

Since we have chosen to consider elements which are centered, the sum (\ref{dev k7 52}) can be restricted to pairs of partitions without singletons. This leaves us with the following pairs $(\pi,\pi')$: $(\{\{1,2\}\},\{\{3,4,5,6,7\}\})$, $(\{\{1,2\}\},\{\{3,4\},\{5,6,7\}\})$ and those which are deduced from the latter by permuting $3,4,5,6,7$.\\

Let us compute $D_{c}(\{\{1,2\}\},\{\{3,4,5,6,7\}\})$. By the second assertion of Lemma \ref{c sigma pi pi}, the permutations $\sigma$ which contribute to this term must have $1,2$ on one hand, and $3,4,5,6,7$ on the other hand, in the same cycle. This can occur if $\sigma$ is either a $7$-cycle or the product of a $2$-cycle and a $5$-cycle.

Let us first compute the contribution of $7$-cycles. The coefficient $C(\sigma,\{1,2\},\{3,4,5,6,7\})$ is not the same for all $7$-cycles $\sigma$. We must distinguish between those in which $1$ and $2$ are consecutive and those in which they are not. There are $2!5!$ $7$-cycles in which $1$ and $2$ are consecutive. For each such cycle $\sigma$, $C(\sigma,\{1,2\},\{3,4,5,6,7\})=1$, thanks to (\ref{eddie.vedder.11.9.08}). In a cycle where $1$ and $2$ are not consecutive, there  may be one, two, three or four elements of $\{3,4,5,6,7\}$ between $1$ and $2$. In each case, there are $5!$ cycles, each contributing a factor $e^{-2t}$, thanks to (\ref{eddie.vedder.11.9.08.bis}).

Let us now compute the contribution of products of a transposition and a $5$-cycle. There are $1!4!$ permutations with two cycles, one which contains $1,2$ and the other $3,4,5,6,7$. Each such permutation contributes a factor $c_{52}$.

Altogether, we have found that 
\begin{equation}\label{d25 1}
D_{c}(\{\{1,2\}\},\{\{3,4,5,6,7\}\})=24(c_{52}+10(1+2e^{-2t})).
\end{equation}

Let us now compute $D_{c}(\{\{1,2\}\},\{\{3,4\},\{5,6,7\}\})$. By the second assertion of Lemma \ref{c sigma pi pi}, there are five possibilities for the partition $\{\{\sigma\}\}$ underlying a permutation $\sigma$ which contributes to this coefficient. We study them one after the other.

${\scriptstyle \bullet} \{\{\sigma\}\}=\{\{1,2,3,4,5,6,7\}\}$. Since, by the third assertion of Lemma \ref{c sigma pi pi}, any two elements of $\{3,4,5,6,7\}$ which are consecutive in $\sigma$ must be in the same block of $\pi'=\{\{3,4\},\{5,6,7\}\}$, no element of $\{3,4\}$ can be consecutive to an element of $\{5,6,7\}$ in $\sigma$. Since there are only two $a$'s, the only possibility is that $3$ and $4$ on one hand, and $5$, $6$, and $7$ on the other hand, are consecutive in $\sigma$ and separated by $1$ and $2$. There are $2!2!3!$ $7$-cycles with this property. Each of them contributes to the sum with a factor $1-e^{-2t}$, according to (\ref{eddie.vedder.11.9.08.bis}).

${\scriptstyle \bullet} \{\{\sigma\}\}=\{\{1,2\},\{3,4,5,6,7\}\}$. By the third assertion of Lemma \ref{c sigma pi pi}, these permutations do not contribute.

${\scriptstyle \bullet} \{\{\sigma\}\}=\{\{1,2,3,4\},\{5,6,7\}\}$. There are two possible structures for the $4$-cycle of $\sigma$ in this case. Either the $a$'s and the $b$'s are consecutive, or they are intertwined. In the first situation, there are $2!2!2!$ permutations, each of which contributes $c_{43}$, thanks to (\ref{eddie.vedder.11.9.08}). In the second situation, there are $2!2!$ permutations, because of a higher symmetry, each of which contributes $e^{-2t}c_{43}$, thanks to (\ref{eddie.vedder.11.9.08.bis}).

${\scriptstyle \bullet} \{\{\sigma\}\}=\{\{1,2,5,6,7\},\{3,4\}\}$. Again, there are two possible structures for the $5$-cycle of $\sigma$, depending on whether the $a$'s are consecutive or not. There are $2!3!$ permutations where they are, and each contributes $c_{52}$. There are also $2!3!$ permutations where they are not, each contributing $e^{-2t} c_{52}$.

${\scriptstyle \bullet} \{\{\sigma\}\}=\{\{1,2\},\{3,4\},\{5,6,7\}\}$. This is the simplest situation. There are $2$ permutations with this cycle structure and each contributes $c_{322}$.

Finally, 
\begin{equation}\label{d52 2}
D_{c}(\{\{1,2\}\},\{\{3,4\},\{5,6,7\}\})=2\left(c_{322}+2(2+e^{-2t}) c_{43} + 6(1+e^{-2t})c_{52} + 12(1-e^{-2t})\right).
\end{equation}

Let us perform the same kind of computations on 
$$k_7(a_1,a_2,a_3,b_1,b_2,b_3,b_4)=\sum_{\pi\in\mc{P}(\{1,2,3\}), \pi'\in \mc{P}(\{4,5,6,7\})} D_{c}(\pi,\pi') \varphi_{\pi}(a_1,a_2,a_3) \varphi_{\pi'}(b_1,b_2,b_3,b_4).$$
Since our variables are centered, the only pairs of partitions which occur in the sum are $(\{\{1,2,3\}\},\{\{4,5,6,7\}\})$, $(\{\{1,2,3\}\},\{\{4,5\},\{6,7\}\})$ and those which are deduced from the latter by permuting $4,5,6,7$.\\

Let us compute $D_{c}(\{\{1,2,3\}\},\{\{4,5,6,7\}\})$. We shall again use formulae \eqref{eddie.vedder.11.9.08}, \eqref{eddie.vedder.11.9.08.bis} and \eqref{eddie.vedder.11.9.08.third} several times, without citing them every time.
The permutations which contribute to this coefficient are $7$-cycles and products of a $3$-cycle and a $4$-cycle. As before, all $7$-cycles do not contribute in the same way. If the $a$'s are consecutive, which is the case for $3!4!$ $7$-cycles, the contribution is simply $1$. If two $a$'s  are consecutive and the third is on its own, the $7$-cycle contributes $e^{-2t}$. In this case, there can be one, two or three $b$'s between the isolated $a$ and the pair of consecutive $a$'s, in the cyclic order. In each case, there are $3!4!$ possible $7$-cycles. Finally, the three $a$'s can be isolated. This happens in $3!4!$ $7$-cycles, and each contributes $e^{-3t}$, thanks to (\ref{eddie.vedder.11.9.08.third}).
So far, we have a contribution of $144(1+3e^{-2t}+e^{-3t})$. The contribution of products of a $3$-cycle and a $4$-cycle is much simpler to compute: it is $2!3!c_{43}$. We find 
\begin{equation}\label{d34 1}
D_{c}(\{\{1,2,3\}\},\{\{4,5,6,7\}\})=12(c_{43}+12(1+3e^{-2t}+e^{-3t})).
\end{equation}

Let us finally compute $D_{c}(\{\{1,2,3\}\},\{\{4,5\},\{6,7\}\})$. Again, by Lemma \ref{c sigma pi pi}, there are five possibilities for the partition $\{\{\sigma\}\}$, which we examine one after the other.

${\scriptstyle \bullet} \{\{\sigma\}\}=\{\{1,2,3,4,5,6,7\}\}$. No element of $\{4,5\}$ can be consecutive with an element of $\{6,7\}$ in $\sigma$. Still, there are several possibilities. Let us first consider the $7$-cycles where $4,5$ on one hand and $6,7$ on the other hand are consecutive. These two groups must be separated by $a$'s. There are $2!2!2!3!$ such $7$-cycles, and each contributes for $1-e^{-2t}$, according to (\ref{eddie.vedder.11.9.08.bis}). Since there are only three $a$'s, one at least of the two pairs $\{4,5\}$ and $\{6,7\}$ must be consecutive. However, it can happen that one is not. This happens in $2!2!2!3!$ $7$-cycles, and according to (\ref{eddie.vedder.11.9.08.third}), each contributes for $e^{-2t}-e^{-3t}$. 

${\scriptstyle \bullet} \{\{\sigma\}\}=\{\{1,2,3\},\{4,5,6,7\}\}$. These permutations do not contribute.

${\scriptstyle \bullet} \{\{\sigma\}\}=\{\{1,2,3,4,5\},\{6,7\}\}$. As usual by now, there are two possibilities in the $5$-cycle of $\sigma$. Either the two $b$'s are consecutive, which happens in $2!3!$ cases with the contribution $c_{52}$, or they are not. This happens in $2!3!$ cases, and each case contributes for $e^{-2t}c_{52}$.

${\scriptstyle \bullet} \{\{\sigma\}\}=\{\{1,2,3,6,7\},\{4,5\}\}$. By symmetry, this contribution is equal to the one above.

${\scriptstyle \bullet} \{\{\sigma\}\}=\{\{1,2,3\},\{4,5\},\{6,7\}\}$. There are $2$ permutations, each contributing for $c_{322}$.

Finally,
\begin{equation}\label{d34 2}
D_{c}(\{\{1,2,3\}\},\{\{4,5\},\{6,7\}\})=2(c_{322}+12(1+e^{-2t})c_{52} + 24 (1-e^{-3t})).
\end{equation}

Let us summarize our results. We have proved that, if there exists a $t$-free cumulant of order $7$, whose associated $7$-linear form is denoted by $k_7$, then for all centered $a_1,a_2,a_3\in \A$ and $b_1,\ldots,b_5 \in \mc{B}$, the following equalities hold:
\begin{align*}
k_7(a_1,a_2,b_1,b_2,b_3,b_4,b_5) &=24(c_{52}+10(1+2e^{-2t})) \varphi(a_{1}a_{2})\varphi(b_{1}b_{2}b_{3}b_{4} b_{5}) \\
& \hskip -3cm + 2\left(c_{322}+2(2+e^{-2t}) c_{43} + 6(1+e^{-2t})c_{52} + 12(1-e^{-2t})\right) \varphi(a_{1}a_{2})\varphi(b_{1}b_{2})\varphi(b_{3}b_{4}b_{5})\\
& \hskip -3cm + 2\left(c_{322}+2(2+e^{-2t}) c_{43} + 6(1+e^{-2t})c_{52} + 12(1-e^{-2t})\right) \varphi(a_{1}a_{2})\varphi(b_{1}b_{3})\varphi(b_{2}b_{4}b_{5})\\
& \hskip -3cm + \ldots,
\end{align*}
where all partitions of $\{b_1,b_2,b_3,b_4,b_5\}$ into a pair and a triple appear, and

\begin{align*}
k_7(a_1,a_2,a_3,b_1,b_2,b_3,b_4) &=12(c_{43}+12(1+3e^{-2t}+e^{-3t})) \varphi(a_{1}a_{2}a_3)\varphi(b_{1}b_{2}b_{3}b_{4}) \\
&\hskip -2cm  + 2(c_{322}+12(1+e^{-2t})c_{52} + 24 (1-e^{-3t})) \varphi(a_{1}a_{2}a_3)\varphi(b_{1}b_{2})\varphi(b_{3}b_{4})\\
&\hskip -2cm  + 2(c_{322}+12(1+e^{-2t})c_{52} + 24 (1-e^{-3t})) \varphi(a_{1}a_{2}a_3)\varphi(b_{1}b_{3})\varphi(b_{2}b_{4})\\
&\hskip -2cm  + 2(c_{322}+12(1+e^{-2t})c_{52} + 24 (1-e^{-3t})) \varphi(a_{1}a_{2}a_3)\varphi(b_{1}b_{4})\varphi(b_{2}b_{3}).
\end{align*}
Since $k_7$ is a $t$-free cumulant, these two expressions must be zero for all choices of $a$'s and $b$'s. %Since the distribution of the $a$'s and $b$'s is unspecified, 
Since the joint  distributions of the $a$'s and of the $b$'s are both arbitrary among those of families of centered elements,
this implies that the coefficients which appear in these equalities in front of the various products of moments of $a$'s and $b$'s must vanish. This implies the following relations:
\begin{align*}
c_{52}&=-10(1+2e^{-2t}),\\
c_{43}&=-12(1+3e^{-2t}+e^{-3t}),\\
c_{322}&=-2(2+e^{-2t}) c_{43} - 6(1+e^{-2t})c_{52} - 12(1-e^{-2t}),\\
c_{322}&=-12(1+e^{-2t})c_{52} - 24 (1-e^{-3t}).
\end{align*}

It does not take a long computation to see that the two expressions of $c_{322}$ are different, since the first involves $e^{-5t}$, whereas the second does not. We leave it to the reader to check that the difference between the two values of $c_{322}$ that we have obtained is equal to $24 e^{-3t} (1-e^{-t})^2$. This quantity vanishes only for $t=0$ or $t=+\infty$.
\end{proof}

In order to prove that $t$-free cumulants of order at most $6$ exist, we are going to construct them. We prove first a lemma which settles the problem of the coefficients $c_\lambda$ for the partitions $\lambda$ whose smallest part is $1$. 

Let us introduce some notation. Let $\mu=(\mu_1\geq \ldots \geq \mu_r)$ be a partition of some integer $n$. We denote by $\ell(\mu)$ the number of non-zero parts of $\mu$ and we write $\mu \vdash \!\! \vdash n$ if $\mu_{\ell(\mu)}\geq 2$, that is, if $\mu$ has no part equal to $1$. Let $i\geq 1$ an integer. We denote by $\mu+\delta_i$ the partition of $n+1$ whose parts are $\mu_1,\ldots,\mu_{i-1},\mu_i+1,\mu_{i+1},\ldots,\mu_r$, rearranged in non-increasing order. If $i>\ell(\mu)$, then $\mu+\delta_i$ is simply the partition $\mu$ to which a part equal to $1$ has been appended. 

\begin{proposition}\label{cond part avec 1} \label{determine} Let $n\geq 2$ be an integer. Choose $t\in [0,+\infty]$. A collection $(c_\lambda)_{\lambda\vdash n}$ is a $t$-free cumulant if and only if the following two conditions hold:\\
1. The relation (\ref{cum 1}) is satisfied for all $m_{1},\ldots,m_{n}$ which are centered.\\
2. For all $\mu\vdash n-1$,  
\begin{equation}\label{part avec 1}
c_{\mu+\delta_{\ell(\mu)+1}}=-\sum_{i=1}^{\ell(\mu)} \mu_i c_{\mu+\delta_i}.
\end{equation}

Moreover, a collection of complex numbers  $(c_\lambda)_{\lambda \vdash\!\!\vdash n}$ which satisfies  (\ref{cum 1}) for all $m_1,\ldots,m_n$ which are centered can be extended in a unique way into a $t$-free cumulant of order $n$.
\end{proposition}

When $\sigma$ is a permutation of $\{1,\ldots,n\}$, let us denote by $[\sigma]$ the partition of the integer $n$ given by the lengths of the cycles of $\sigma$. 

\begin{proof} Let $c$ be a $t$-free cumulant of order $n$. In order to check that (\ref{part avec 1}) is satisfied, let us choose $m_1,\ldots,m_{n-1}$ in some probability space $(\mc{M},\varphi)$ and write the fact that $k_n(m_1,\ldots,m_{n-1},1)=0$. We find
\begin{equation}\label{mu 1}
\sum_{\lambda\vdash n}c_\lambda  \sum_{\substack{\sigma \in \Sym_n\\ [\sigma]=\lambda}} \varphi_\sigma(m_1,\ldots,m_{n-1},1)=0.
\end{equation}
Let $r_n:\Sym_n\to \Sym_{n-1}$ denote the following function: for all $\sigma\in \Sym_n$, $r_n(\sigma)$ is the permutation of $\{1,\ldots,n-1\}$ obtained by removing $n$ from the cycle of $\sigma$ which contains it. 
For each $\sigma$, we have the equality $\varphi_\sigma(m_1,\ldots,m_{n-1},1)=\varphi_{r_n(\sigma)}(m_1,\ldots,m_{n-1})$. Now a permutation $\tau\in \Sym_{n-1}$ has exactly $n$ preimages by $r_n$. Moreover, if $[\tau]=\mu=(\mu_1\geq \ldots \geq \mu_{\ell(\mu)}>0) \vdash n-1$, then all preimages of $\tau$ belong to one of the conjugacy classes $\mu+\delta_i$ for $i=1,\ldots,\ell(\mu)+1$. Finally, $r_n^{-1}(\tau)$ contains exactly one element of $\mu+\delta_{\ell(\mu)+1}$ and $\mu_i$ elements of $\mu+\delta_i$ for $i=1,\ldots,\ell(\mu)$. We can thus rewrite (\ref{mu 1}) as follows:
\begin{equation}\label{mu 2}
\sum_{\mu\vdash n-1} \left(\sum_{i=1}^{\ell(\mu)}  \mu_i c_{\mu+\delta_i} + c_{\mu+\delta_{\ell(\mu)+1}}  \right)  \sum_{\substack{\tau \in \Sym_{n-1}\\ [\tau]=\mu}} \varphi_\tau(m_1,\ldots,m_{n-1})=0.
\end{equation}
Since the distribution of $(m_1,\ldots,m_{n-1})$ is arbitrary, all the coefficients between the brackets must vanish. It follows that (\ref{part avec 1}) is satisfied.

Conversely, let $(c(\sigma))_{\sigma\in \Sym_{n}}$ be a collection which satisfies (\ref{cum 1}) for centered elements and (\ref{part avec 1}). Then, by the computation that we have just done, this collection satisfies (\ref{cum 2}) and hence, by multilinearity, (\ref{cum 1}) for arbitrary elements.

Let us prove the last assertion. For any $\lambda\vdash n$ with at least one part equal to $1$, the relation (\ref{part avec 1}) expresses the value of $c_\lambda$ in terms of $c_{\lambda'}$ for partitions $\lambda'$ of $n$ which have strictly less parts equal to $1$ than $\lambda$. The collection $(c_\lambda)_{\lambda\vdash n}$ is thus completely and uniquely determined by $(c_\lambda)_{\lambda \vdash\!\!\vdash n}$. The fact that the resulting collection is a $t$-free cumulant is granted by the first part of the proposition.
\end{proof}

The last result simplifies greatly the search for $t$-free cumulants, since it allows one to restrict to centered elements and partitions in parts at least equal to $2$. We apply it to find cumulants of order less than $6$.

\begin{proof}(Theorem \ref{existence}) Let us prove that there exist $t$-free cumulants up to order $6$. We proceed by first establishing enough conditions that their coefficients must satisfy, in order to determine these coefficients. Then, we check that we actually have a $t$-free cumulant.

We will always normalize our cumulants by the condition $c_n=1$.

${\scriptstyle\bullet} \; n=2.$ By Proposition \ref{determine}, the condition $c_2=1$ suffices to determine the whole cumulant, and $c_{11}=-1$. The relation (\ref{eddie.vedder.11.9.08}) implies that we have indeed got a $t$-free cumulant.

${\scriptstyle\bullet} \; n=3.$ Again, the condition $c_3=1$ determines completely the cumulant. Using (\ref{part avec 1}), we find $c_{21}=-2$ and $c_{111}=4$. The relation (\ref{eddie.vedder.11.9.08}) implies again that the collection thus obtained is a $t$-free cumulant. Indeed, by \eqref{eddie.vedder.11.9.08}, the product of any three centered elements, one being $t$-free with the two others, is centered. Hence, our collection satisfies (\ref{cum 1}) on centered elements.

${\scriptstyle\bullet} \; n=4.$ This is the first case where the relation $c_4=1$ does not suffice determine the cumulant. Indeed, we must compute $c_{22}$. For this, let us choose in some probability space elements $a_1,a_2,\ldots$ and $b_1,b_2,\ldots$, such that $\{a_{1},a_{2},\ldots\}$ and $\{b_{1},b_{2},\ldots\}$ are $t$-free. We will use this notation again in this proof without redefining it. Let us assume that a $t$-free cumulant $c$ of order $4$ is given and let us compute $D_{c}(\{\{1,2\}\},\{\{3,4\}\})$ (see (\ref{def D})). Again, we shall use formulae \eqref{eddie.vedder.11.9.08} and  \eqref{eddie.vedder.11.9.08.bis}  several times, without citing them every time. There are $4$-cycles which contribute to this coefficient. In $2!2!$ of them, $1$ and $2$ are consecutive and they contribute for $1$ each. In $2!$ others, $1$ and $2$ are not consecutive and each such cycle contributes for $e^{-2t}$. There is also one product of two $2$-cycles, which contributes for $c_{22}$. Finally, $D_{c}(\{\{1,2\}\},\{\{3,4\}\})=c_{22}+2(2+e^{-2t})$. The nullity of this coefficient implies $c_{22}=-2(2+e^{-2t})$. Using (\ref{part avec 1}), we determine the remaining coefficients, and find
\begin{align*}
c_4=1,  \; c_{31}=-3,\; c_{22}=-2(2+e^{-2t}), \; c_{211}=2(5+e^{-2t}), c_{1111}=-6(5+e^{-2t}).
\end{align*}
Now let us check that the collection thus defined satisfies (\ref{cum 1}) for centered elements. Set $k_{4}=\sum_{\sigma} c(\sigma)\varphi_{\sigma}$. If we expand $k_{4}(a_{1},b_{1},b_{2},b_{3})$ according to (\ref{expand}), then all terms involve $\varphi(a_{1})$ and vanish. Now $k_{4}(a_{1},a_{2},b_{1},b_{2})$ also vanishes, because this is how we have chosen the value of $c_{22}$. Finally, we do have a $t$-free cumulant of order $4$. 

${\scriptstyle\bullet} \; n=5.$ Let $c$ be a $t$-free cumulant of order $5$. Let us compute $c_{32}$ by writing the nullity of $D_{c}(\{\{1,2,3\}\},\{\{4,5\}\})$ and using formulae \eqref{eddie.vedder.11.9.08} and \eqref{eddie.vedder.11.9.08.bis}. There are $3!2!$ $5$-cycles in which $4$ and $5$ are consecutive, and they contribute for $1$ each. There are also $3!2!$ $5$-cycles in which they are not consecutive, and each contributes for $e^{-2t}$. There are finally $2!$ products of a $3$-cycle and a $2$-cycle, which contribute for $c_{32}$ each. Hence, we must have $c_{32}=-6(1+e^{-2t})$. Using as usual (\ref{part avec 1}), we find that the other values of $c$ must be
\begin{align*}
&c_5=1, \; c_{41}=-4, \; c_{32}=-6(1+e^{-2t}), \; c_{311}=6(3+e^{-2t}), c_{221}=12(1+e^{-2t}), \\
& c_{2111}=-12(5+2e^{-2t}), \; c_{11111}=48(5+2e^{-2t}).
\end{align*}
The fact that $k_{5}\sum_{\sigma} c(\sigma)\varphi_{\sigma}$ is a cumulant is checked just as in the case $n=4$: the identity $k_{5}(a_{1},b_{1},b_{2},b_{3},b_{4})=0$ is granted by (\ref{expand}) and $k_{5}(a_{1},a_{2},b_{1},b_{2},b_{3})=0$ by the choice of $c_{32}$. 

${\scriptstyle\bullet} \; n=6.$ Let $c$ be a $t$-free cumulant of order $5$. The value of $c_{42}$, deduced as usual from the nullity of $D_{c}(\{\{1,2,3,4\}\},\{\{5,6\}\})$, is $c_{42}=-4(2+3e^{-2t})$. Similarly, $D_{c}(\{\{1,2,3\}\},\{\{4,5,6\}\})=0$ gives us $c_{33}=-3(3 + 6 e^{-2t} +  e^{-3t})$. Finally, $D_{c}(\{\{1,2\}\},\{\{3,4\},\{5,6\}\})=0$ implies $c_{222}=8(7+17e^{-2t}+6e^{-4t})$. The other coefficients follow as usual from (\ref{part avec 1}) and we find
\begin{align*}
&c_{6}=1, \; c_{51}=-5, \; c_{42}=-4(2+3e^{-2t}), \; c_{411}=4(7+3e^{-2t}), \; c_{33}=-3(3+6e^{-2t}+e^{-3t})\\
&c_{321}=6(7+e^{-3t}+12e^{-2t}), \; c_{3111}=-12(14+15e^{-2t}+e^{-3t}), \;c_{222}=8(7+17e^{-2t}+6e^{-4t})\\
&c_{2211}=-8(28+53e^{-2t}+3e^{-3t}+6e^{-4t}), \; c_{21111}=48(21+34e^{-2t}+2e^{-3t}+3e^{-4t})\\
&c_{111111}=-240(21+34e^{-2t}+2e^{-3t}+3e^{-4t}).  
\end{align*}
Let us set $k_{6}=\sum_{\sigma} c(\sigma) \varphi_{\sigma}$. The nullity of $k_{6}(a_{1},b_{1},\ldots,b_{5})$ follows as usual from (\ref{expand}). That of $k_{6}(a_{1},a_{2},b_{1},b_{2},b_{3},b_{4})$ results from the choices of $c_{42}$ and $c_{222}$. Finally, $k_{6}(a_{1},a_{2},a_{3},b_{1},b_{2},b_{3})=0$ is granted by the choice of $c_{33}$.

Nowhere there has been any freedom in the definition of the cumulants. This shows the uniqueness of conjugation-invariant $t$-free cumulants of order at most than $6$. 
\end{proof}

\end{document}